\documentclass[11pt]{article}

\usepackage{amsmath,amsfonts,amssymb} 

\usepackage{amsthm} 

\usepackage{graphicx,caption,subcaption} 
\usepackage{fullpage}

\allowdisplaybreaks 

\begingroup
    \makeatletter
    \@for\theoremstyle:=definition,remark,plain\do{%
        \expandafter\g@addto@macro\csname th@\theoremstyle\endcsname{%
            \addtolength\thm@preskip\parskip
            }%
        }
\endgroup

\newtheorem{thm}{Theorem}[section]
\newtheorem{cor}[thm]{Corollary}

\newtheorem{lemma}[thm]{Lemma}

\newtheorem{definition}[thm]{Definition}

\newtheorem{claim}[thm]{Claim}

\begin{document}

\title{\vspace{-0.5in} Minimum degree ensuring that a hypergraph is hamiltonian-connected}

\author{
{{Alexandr Kostochka}}\thanks{
\footnotesize {University of Illinois at Urbana--Champaign, Urbana, IL 61801
 and Sobolev Institute of Mathematics, Novosibirsk 630090, Russia. E-mail: \texttt {kostochk@math.uiuc.edu}.
 Research 
is supported in part by  NSF  Grant DMS-2153507.}}
\and
{{Ruth Luo}}\thanks{
\footnotesize {University of South Carolina, Columbia, SC 29208. E-mail: \texttt {ruthluo@sc.edu
}.
 Research 
is supported in part by NSF grant DMS-1902808.
}}
\and{{Grace McCourt}}\thanks{University of Illinois at Urbana--Champaign, Urbana, IL 61801, USA. E-mail: {\tt mccourt4@illinois.edu}. Research 
is supported in part by NSF RTG grant DMS-1937241.}}

\date{ \today}
\maketitle

\vspace{-0.3in}

\begin{abstract}
A hypergraph $H$ is {\em hamiltonian-connected} if for any distinct  vertices $x$ and $y$, $H$ contains a hamiltonian Berge path from $x$ to $y$. We find for all $3\leq r<n$, exact lower bounds on minimum degree $\delta(n,r)$ of an $n$-vertex $r$-uniform hypergraph $H$ guaranteeing that $H$ is 
 hamiltonian-connected. It turns out that for $3\leq n/2<r<n$, $\delta(n,r)$ is $1$ less than the degree bound guaranteeing the existence  of a hamiltonian Berge cycle.
Moreover, unlike  for graphs,  for each $r \geq 3$ there exists an $r$-uniform hypergraph that is hamiltonian-connected but does not contain a hamiltonian Berge cycle.

\end{abstract}

\section{Introduction and results}

A hypergraph $H$ is a family of subsets of a ground set. We refer to these subsets as the {\em edges} of $H$ and the elements of the ground set as the {\em vertices} of $H$. We use $E(H)$ and $V(H)$ to denote the set of edges and the set of vertices of $H$ respectively. We say that $H$ is {\em $r$-uniform}  (or an $r$-graph, for short) if every edge of $H$ contains exactly $r$ vertices. A {\em graph} is a 2-graph. The {\em degree} $d_H(v)$ of a vertex $v$ in a hypergraph $H$ is the number of edges containing $v$. When it is clear from the context, we may simply write $d(v)$ to mean $d_H(v)$. The {\em minimum degree}, $\delta(H)$, is the minimum over degrees of all vertices of $H$. 

A {\em hamiltonian cycle} (path) in a graph is a cycle (path) that visits every vertex. A graph is {\em hamiltonian} if it contains a hamiltonian cycle.  Furthermore,  a graph is {\em hamiltonian-connected} if there exists a hamiltonian path between every pair of vertices. 

It is well known that determining whether a  graph is hamiltonian is an NP-complete problem. Sufficient conditions for existence of hamiltonian cycles in graphs have been well-studied. In particular,  the famous  Dirac's Theorem~\cite{dirac} says  that for any $n\geq 3$ each $n$-vertex graph $G$ with $\delta(G)\geq n/2$  contains a hamiltonian cycle. 

Every hamiltonian-connected graph is also hamiltonian, but the converse is not true. 
For example for even $n\geq 4$,  the  complete bipartite graph $K_{n/2,n/2}$ is hamiltonian but not 
hamiltonian-connected.
The  example of $K_{n/2,n/2}$ also shows that
for even $n$, condition $\delta(G)\geq n/2$ does not provide that $G$    is hamiltonian-connected.
 On the other hand, Ore~\cite{ore2} proved that a slightly stronger restriction on minimum degree of a graph implies hamiltonian-connectedness:

\begin{thm}[Ore~\cite{ore2}]\label{o2} Let $n\geq 3$ and  $G$ be an $n$-vertex graph. If $d(u) + d(v) \geq n + 1$ for 
every  $u,v\in V(G)$ with $uv \notin E(G)$, then $G$ is hamiltonian-connected. In particular, if $\delta(G)\geq (n+1)/2$, then $G$ is hamiltonian-connected.
\end{thm}

Note that for odd $n$, the restriction on minimum degree is the same as in Dirac's Theorem.


 

Dirac's Theorem and Theorem~\ref{o2} have been generalized and refined in several directions by Posa~\cite{Posa}, Lick~\cite{Lick} and 
many others. Among generalizations, there were different extensions of the theorems to cycles and paths in hypergraphs, in particular, in
$r$-graphs.


\begin{definition}
A {\bf Berge cycle of} length $s$ in a hypergraph is a list of $s$ distinct vertices and $s$ distinct edges $v_1, e_1, v_2, \ldots,e_{s-1}, v_s, e_s, v_1$ such that $\{v_i, v_{i+1}\} \subseteq e_i$ for all $1\leq i \leq s$ (we always take indices of cycles of length $s$ modulo $s$). We call vertices $v_1, \dots, v_s$ the {\bf defining vertices} of $C$ and write $V(C)=\{v_1, \ldots, v_s\}$, $E(C) = \{e_1, \ldots, e_s\}.$
Similarly, a {\bf Berge path} of length $\ell$ is a list of $\ell+1$ distinct vertices and $\ell$ distinct edges $v_1, e_1, v_2, \ldots, e_{\ell}, v_{\ell + 1}$ such that $\{v_i, v_{i+1}\} \subseteq e_i$ for all $1\leq i \leq \ell$, with {\bf defining vertices} $V(P) = \{v_1, \dots, v_{\ell+1}\}$ and $E(P) = \{e_1, \dots, e_{\ell}\}$. 
\end{definition}

 For simplicity, we will say a hypergraph is {\em hamiltonian} if it contains a hamiltonian Berge cycle, and is {\em hamiltonian-connected} if it contains a hamiltonian Berge path between any pair of vertices.

Approximate bounds on the minimum degree of an $n$-vertex $r$-graph $H$ that provide that $H$ is hamiltonian were obtained for $r\leq \frac{n-4}{2}$ by Bermond, Germa, Heydemann, and Sotteau~\cite{BGHS}; Clemens, Ehrenm\"uller, and Person~\cite{CEP}; and Ma, Hou, and Gao~\cite{MHG}. Coulson and Perarnau~\cite{CP} gave exact bounds in the case  $r = o(\sqrt{n})$ (and large $n$).
The present authors resolved the problem for all $3\leq r<n$:

\begin{thm}[\cite{KLM}]\label{oldmain} Let $n > r \geq 3$. Suppose $H$ is an $n$-vertex, $r$-graph such that (1) $r \leq (n-1)/2$ and $\delta(H) \geq {\lfloor (n-1)/2 \rfloor \choose r-1} + 1$, or (2) $r \geq n/2$ and $\delta(H) \geq r$.
Then $H$ contains a hamiltonian Berge cycle. 
\end{thm}
The inequalities in this result are best possible for all $3\leq r<n$. Very recently, Salia~\cite{SN} 
proved sharp results of P\' osa type for Berge hamiltonian cycles.
He 
 described the sequences $(d_1,\ldots,d_n)$ with $d_1\leq d_2\leq\ldots\leq d_n$ of two types: 
(a) for  $r<n/2$ every $n$-vertex $r$-graph with degree sequence $(d'_1,\ldots,d'_n)$ such that $d'_i>d_i$ for all $i$ has a hamiltonian Berge cycle and also  (b)  every $n$-vertex hypergraph with degree sequence $(d'_1,\ldots,d'_n)$ such that $d'_i>d_i$ for all $i$ has a hamiltonian Berge cycle. The first of these nice results implies Part (a) of Theorem~\ref{oldmain} for odd $n$.

Since we consider mostly Berge cycles and paths, from now on, we will drop the word ``Berge" and simply use {\em cycle} and {\em path} to refer to a Berge cycle and a Berge path, respectively.

Note that while every hamiltonian-connected graph is hamiltonian, this is not true for $r$-graphs when $3\leq r<n$.
In the next section, for every $3\leq r<n$ we present a hamiltonian-connected $r$-graph that has no hamiltonian cycles.

The main result of  this paper is the following. 

\begin{thm}\label{main} Let $n \geq r \geq 3$. Suppose $H$ is an $n$-vertex $r$-graph such that\\ (1) $r \leq n/2$ and $\delta(H) \geq {\lfloor n/2 \rfloor \choose r-1} + 1$, $\quad$ or (2) $n-1\geq r > n/2\geq 3$  and $\delta(H) \geq r-1$,\\ or (3) $r=3$, $n=5$ and $\delta(H) \geq 3$.
Then $H$ is hamiltonian-connected. 
\end{thm}

Note that the conditions in Theorem~\ref{main} for $3\leq r\leq n/2$ and even $n$ are stronger than 
in Theorem~\ref{oldmain}, for $3\leq r\leq n/2$ and odd $n$  are the same, and for $3\leq n/2< r\leq n-1$  are weaker than 
in Theorem~\ref{oldmain}. These bounds are sharp, and extremal examples will be given in the next section.

Similarly to~\cite{KLM}, we elaborate the idea of Dirac~\cite{dirac} of choosing a longest cycle plus a longest path. We also use a series of lemmas on subsets of edges and vertices in graph paths.

The structure of the paper is as follows. In the next section, we show extremal examples for Theorems~\ref{oldmain}
and~\ref{main} and also examples of hamiltonian-connected $r$-graphs that have no  hamiltonian cycles.
In Section~3 we prove lemmas on subsets of graph paths. In Section~4 we set up the main proofs for all cases:
 we define ``best" extremal substructures in possible counter-examples to our theorem
 and prove some properties of such substructures. In the subsequent three sections, we analyze all possible cases that can arise in counter-examples, and settle these cases. We finish the paper with some concluding remarks.

\section{Examples}

\subsection{Examples for Theorems~\ref{oldmain} and~\ref{main} }

For all $n>3$ and $3\leq r \leq (n-1)/2$, let $H_1=H_1(n,r)$ be the $r$-graph 
formed by a clique $Q$ of size $\lceil \frac{n+1}{2}\rceil$ and a clique $R$ of size $\lfloor \frac{n+1}{2}\rfloor$
 sharing exactly one vertex.
Then $\delta(H_1) = { \lfloor \frac{n-1}{2}\rfloor \choose r-1}$, and $H_1$ is non-hamiltonian because it has  a vertex
whose deletion disconnects the $r$-graph.

Another example for $3\leq r \leq (n-1)/2$, is the $r$-graph $H_2=H_2(n,r)$ whose vertex set is
$A\cup B$ where $|A|=\lceil \frac{n+1}{2}\rceil$, $|B|=\lfloor \frac{n-1}{2}\rfloor$, $A\cap B=\emptyset$ and whose edges are sets
$X\subset A\cup B$ with $|X|=r$ and $|X\cap A|\leq 1$. Again, $\delta(H_2) = { \lfloor \frac{n-1}{2}\rfloor \choose r-1}$. Also, each cycle in $H_2$ has no two consecutive vertices in $A$. Since $|A|>n/2$, this yields that $H_2$ is not hamiltonian.

For $n/2\leq r \leq n-1$, $H_3=H_3(n,r)$ is obtained by removing a single edge from any $r$-regular $r$-graph. Then $\delta(H_3) = r-1$ and $H_3$ has $n-1$ edges. Hence $H_3$ cannot have a hamiltonian  cycle. 

The $r$-graphs above show sharpness of the bounds in Theorem~\ref{oldmain}. The following slight modifications of them show sharpness of the bounds in Theorem~\ref{main}.

For all $n>3$ and $3\leq r \leq n/2$, let $H'_1=H'_1(n,r)$ be the $r$-graph 
formed by a clique $Q$ of size $\lceil \frac{n+2}{2}\rceil$ and a clique $R$ of size $\lfloor \frac{n+2}{2}\rfloor$
 sharing exactly two vertices, say $x$ and $y$. Then $\delta(H'_1) = { \lfloor \frac{n}{2}\rfloor \choose r-1}$, and $H'_1$
has no hamiltonian $x,y$-path, since any $x,y$-path should miss either $Q-\{x,y\}$ or $R-\{x,y\}$.

Another example for $3\leq r \leq n/2$, is the $r$-graph $H'_2=H'_2(n,r)$ whose vertex set is
$A\cup B$ where $|A|=\lceil \frac{n}{2} \rceil$ and $|B|=\lfloor\frac{n}{2}\rfloor$, $A\cap B=\emptyset$ and whose edges are sets
$X\subset A\cup B$ with $|X|=r$ and $|X\cap A|\leq 1$. Now $\delta(H'_2) = { \lfloor \frac{n}{2} \rfloor \choose r-1}$. Also,
for distinct $x,y\in B$ each $x,y$-path in $H'_2$ has no two consecutive vertices in $A$. Since $|A|\geq n/2$, this yields that $H'_2$ has no hamiltonian $x,y$-path.

For $r > n/2$, let $H'_3=H'_3(n,r)$ be obtained from $H_3(n,r)$ by removing any edge.
 Then $\delta(H'_3) = r-2$ and $H'_3$ has $n-2$ edges. Hence $H'_3$ cannot have any hamiltonian  path. 
 
 For $r=3$, $n=5$, let $V(H_4) = \{1,2,3,4,5\}$ and $E(H_4) =\{ \{1,5,2\}, \{1,5,3\}, \{1,5,4\}, \{2,3,4\}\}$. Then $\delta(H_4) = 2$ but there is no hamiltonian path from $1$ to $5$.

\subsection{Hamiltonian-connected $r$-graphs with no hamiltionian cycles}\label{hcnh}

By the $(n,r)$-{\em tight cycle } $C(n,r)$ we denote the $r$-graph with vertex set $V=\{v_1,\ldots,v_n\}$ and edge  
set $E=\{e_1,\ldots,e_n\}$, where $e_i=\{v_i,v_{i+1},\ldots,v_{i+r-1}\}$ for all $i=1,\ldots,n$ and indices count modulo $n$.

Our example $C'(n,r)$ is obtained from $C(n,r)$ by deleting one edge. Since $C'(n,r)$ has $n-1$ edges, it has 
no hamiltionian cycles. We claim that for $3\leq r<n$, $C'(n,r)$ is hamiltonian-connected.

Indeed, by symmetry we may assume that 
  we need a hamiltionian $v_1,v_h$-path and that
  we have deleted $e_j$ from 
$C(n,r)$. Also by symmetry, we may assume that $h\leq j+1\leq n$. We construct a hamiltionian $v_1,v_h$-path slightly differently for odd $h$, for even $h\geq 4$ and for $h=2$.
In all cases,  the subpath from $v_n$ to $v_h$ will be
\[P_2 = v_n,e_{n-1},v_{n-1},e_{n-2},
\ldots,
\ldots,v_{j+2},e_{j+1},v_{j+1},e_{j-1},v_{j},e_{j-2},v_{j-1},\ldots, e_{h},v_{h+1},e_{h-1},v_h.\]

Our final hamiltonian $v_1, v_h$-path will be of the form $P_1 \cup P_2$
(see Figure~\ref{fig:hc})
 where the  subpath $P_1$ is as follows:

\begin{figure}
    \centering
    \includegraphics[height=1.9in]{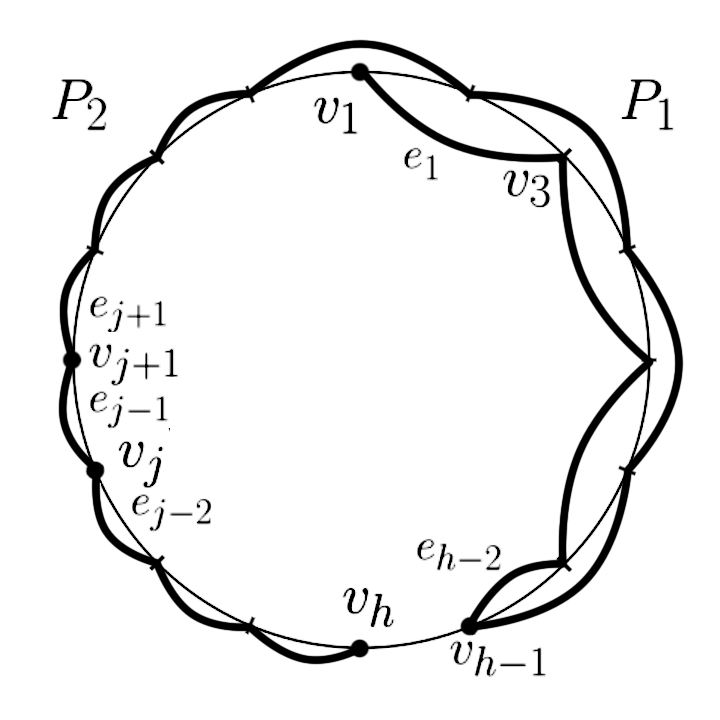} \qquad 
    \includegraphics[height=1.9in]{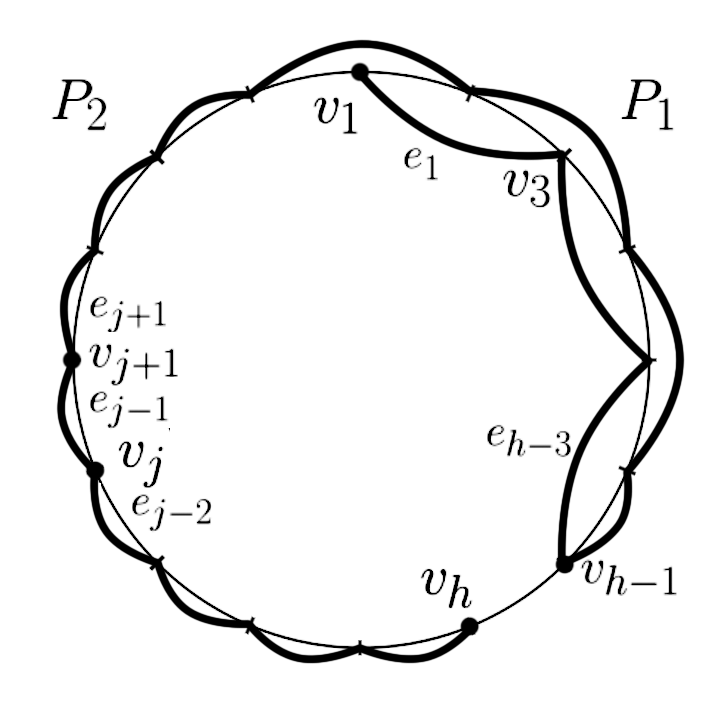} \qquad
     \includegraphics[height=1.9in]{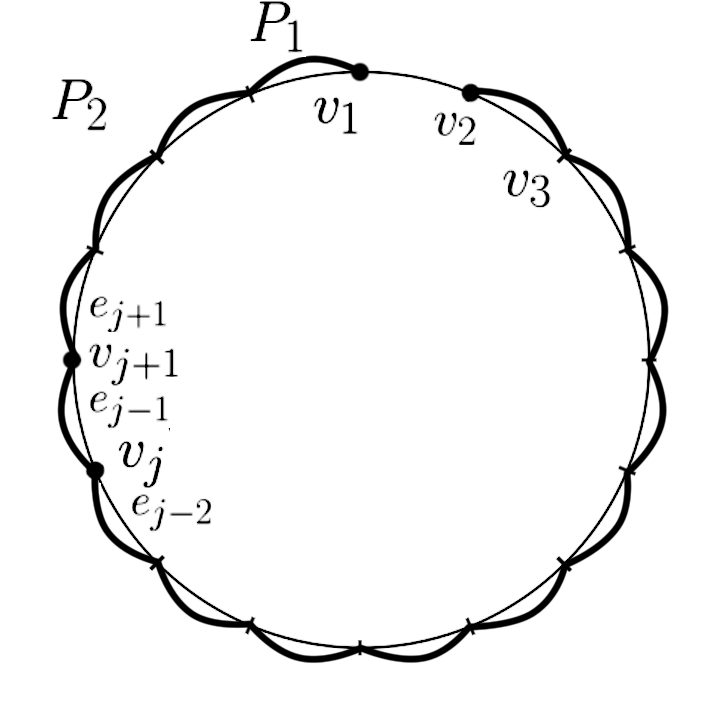}

      \caption{The three cases of a hamiltonian $v_1,v_h$-path in $C'(n,r)$.}
    \label{fig:hc}
\end{figure}

If $h$ is odd, then
\[P_1 = v_1,e_1,v_3,e_3,v_5,\ldots,v_{h-2},e_{h-2},v_{h-1},e_{h-3},v_{h-3},e_{h-5},\ldots,e_2,v_2, e_n, v_n.\] If $h$ is even and $h>2$, then \[P_1 = v_1,e_1,v_3,e_3,v_5,\ldots,v_{h-1},e_{h-2},v_{h-2},e_{h-4},v_{h-4},e_{h-6},\ldots,e_2,v_2, e_n, v_n.\] Finally if $h=2$, then $P_1 = v_1, e_n, v_n$.

\section{Lemmas  on graph paths}

In this section we derive some properties of subsets of graph paths that  will be heavily used in our proofs.
The reader can skip their proofs at the first reading.

\begin{lemma}\label{verc2} Let $Q=v_1,e_1,\ldots,e_{s-1},v_s$ be a graph path. Let $A$ and $B$ be nonempty subsets of $V(Q)$ such that $A$ is an independent set, $B - A \neq \emptyset$, and for each $v_i\in A$ and $v_j\in B-A$, $|i - j| \geq q\geq 1$. Then 

(i) If $q\geq 2$, then $s \geq 2|A| + |B-A| + q-2$ (and therefore $|A| \leq (s-|B-A|-q+2)/2$).

(ii) If $q = 1$, then $s \geq 2|A| + |B-A| -2$.

Moreover, if $B$ is also an independent set, then $s\geq 2|A| + 2|B-A| + q - 3$. 
\end{lemma}

\begin{proof} 
%
%
%
%
%
 
Let $v_j \in B-A$. Without loss of generality, we may suppose there exists a vertex $v_i \in A$ such that $i < j$ and $v_k \notin A \cup B$ for all $i < k <j$. Then $V_1:=\{v_{i+1}, \ldots, v_{j-1}\}$ is a set of of at least $q-1$ vertices which does intersect $A \cup B$. Similarly, if there exists $v_{i'} \in A$ such $i < j < i'$ (and $v_k \notin A \cup B$ for all $j<k<i'$), then $V_2:=\{v_{j+1}, \ldots, v_{i'-1}\}$ also contains at least $q-1$ vertices and does not intersect $A \cup B$. In this case, set $V'  = V_1 \cup V_2$. Otherwise, set $V' = V_1$. 

For each $v_{k} \in A - \{v_s\}$, $v_{k+1}$ does not intersect $A \cup B$, and only one $v_{k+1}$, namely $v_{i+1}$, is in $V'$. Therefore
 \[s \geq |A| + |B - A| + |\{v_{k+1}: v_k \in A, k \notin \{i, s\}\}| + |V'|.\]
 If $V' = V_1 \cup V_2$, then $s \geq  2|A| + |B-A| -2 + 2q -2$ which is at least $2|A| + |B-A| + q-2$ if $q \geq 2$, and at least $2|A| + |B-A| + 2$ if $q = 1$.
 
 If $V' = V_1$, then in this case $v_s \notin A$, so we have $s \geq |A| + |B-A| + (|A|-1) + q-1 = 2|A| + |B-A| +q-2$.

 Suppose now that $B$ is also an independent set, and let $v_j \in B-A$. Again we may suppose there exists $v_i \in A$ with $i < j$. Between $v_j$ and $v_i$ there is a set $V'$ of at least $q-1$ vertices not in $A \cup B$, and for any $v_k \in A \cup B$, $v_{k+1} \notin A \cup B$. Therefore
 
 \[s \geq |A \cup B| + |\{v_{k+1}: v_k \in A \cup B, k \notin \{i, s\}\}| + |V'|\]\[\geq |A| + |B-A| + (|A| + |B-A| - 2) + q-1 = 2|A| + 2|B-A| + q - 3.\] 
 \end{proof}

\begin{lemma}\label{indep2} Let $q\geq 2$ and $s>a\geq 1$. Let 
 $Q=v_1,e_1,\ldots,e_{s-1},v_s$ be a graph path, and $I$ be a non-empty independent subset of
$\{v_1, \ldots, v_s\}$. If  $A'$ is a set of $a$  edges of $Q$ such that the distance in $Q$ 
 from any edge in $A'$ to any vertex in $I$ is at least $q$, then $|I|\leq  \lfloor\frac{s - a-q+1}{2} \rfloor$ if $q \geq 2$, and $|I|\leq  \lfloor \frac{s-a+1}{2} \rfloor = \lceil \frac{s - a}{2} \rceil$ if $q = 1$.
\end{lemma}

\begin{proof} 
Applying Lemma~\ref{verc2} with $A = I$ and $B = \bigcup_{\{i:e_i \in A'\}}\{v_i,v_{i+1} \}$ (so $|B| \geq a+1$) gives the desired bounds.
 \end{proof}

\begin{lemma}\label{ver-new} Let $Q=v_1,e_1,\ldots,e_{s-1},v_s$ be a graph path. Let $A$ and $B$ be nonempty subsets in $V(Q)$ such that 
\begin{equation}\label{dis}
\mbox{ for each $v_i\in A$ and $v_j\in B$, either $i=j$ or $|i - j| \geq q\geq 2$. }
\end{equation}

(i)  If $A=B$, then $s\geq 1+q(|A|-1)$ with equality only if $A=\{v_1,v_{1+q},v_{1+2q},\ldots,v_s\}$.

(ii) If $B\neq A$, then $s \geq |A| + |B| + q-2$ with equality only if $A\subset B$ or $B\subset A$.
\end{lemma}

\begin{proof} Part (i) is obvious. We prove (ii) by induction on $|A\cap B|$.

If $A\cap B=\emptyset$, then $Q$ contains $|A|+|B|$ vertices in $A\cup B$ and at least $q-1$ vertices outside of 
$A\cup B$ between $A$ and  a closest to $A$ vertex in $B$.

Suppose now that (ii) holds for all $A'$ and $B'$ with $|A'\cap B'|<t$ and that $|A\cap B|=t$, say $v_i\in A\cap B$.
By symmetry, we may assume $|A|\leq |B|$.  If $A=\{v_i\}$, then $Q$ has $|B|-1$ vertices in $B-A$ and at least $q-1$
vertices between $v_i$ and a closest to $v_i$ vertex in $B-A$ (such a vertex exists since $B\neq A$). Thus,
$s\geq 1+|B-A|+q-1=|B|+|A|-2$, as claimed.

Finally, suppose $|A|\geq 2$. By definition, $(A\cup B)\cap \{v_{i-q+1}, v_{i-q+2},\ldots,v_{i+q-1}\}=\{v_i\}$.
So since $|B|\geq |A|\geq 2$, $i\geq q+1$ or $i\leq s-q$ (or both). By symmetry, assume $i\geq q+1$.
Define $e'_i=v_{i-q+1}v_{i+1}$ and let $A'=A-v_i$, $B'=B-v_i$, and 
$Q'=v_1,e_1,\ldots,v_{i-q+1},e'_i,v_{i+1},e_{i+1},\ldots,e_{s-1},v_s$. By definition, $A'$ and $B'$ satisfy~\eqref{dis}.
So, by induction, $|V(Q')|\geq |A'|+|B'|+q-2$ with equality only if $A'\subset B'$. Hence
$$s\geq q+|V(Q')|\geq q+ (|A|-1)+(|B|-1)+q-2=|A|+|B|+2(q-2),$$
 with equality only if $A\subset B$. Since $q\geq 2$, this proves (ii).
\end{proof}

\begin{lemma}\label{ed-new} Let $Q=v_1,e_1,\ldots,e_{s-1},v_s$ be a graph path. Let $A'$ and $B'$ be nonempty subsets of $E(Q)$ such that 
\begin{equation}\label{dise}
\mbox{ for each $e_i\in A'$ and $e_j\in B'$, either $i=j$ or $|i - j| \geq q\geq 2$. }
\end{equation}

(i)  If $A'=B'$, then $s-1\geq 1+q(|A'|-1)$ with equality only if $A'=\{e_1,e_{1+q},e_{1+2q},\ldots,e_{s-1}\}$.

(ii) If $B'\neq A'$, then $s-1 \geq |A'| + |B'| + q-2$ with equality only if $A'\subset B'$ or $B'\subset A'$.
\end{lemma}

\begin{proof}
Let $A = \{v_i : e_i \in A'\}$ and $B = \{v_i : e_i \in B'\}$. Since $v_s \notin A\cup B$, the sets $A, B$ are vertex subsets of  the path $Q' = v_1, e_1, \dots, e_{s-2}, v_{s-1}$. So,  Lemma~\ref{ver-new} applied to $A,B$ and $Q'$ yields 
 the desired bounds.
\end{proof}

\begin{lemma}\label{consecpath2}Let $Q=v_1,e_1, \ldots,e_{s-1}, v_{s}$ be a graph path. Suppose $F\subset E(Q)$ and $f=|F|$. Let $A$ and $B$ be subsets of $\{v_1, \ldots, v_{s}\}$ that are vertex-disjoint from all edges in $F$ and such that
\begin{equation}\label{dis2}
\mbox{ for each $v_i\in A$ and $v_j\in B$, either $i=j$ or $|i - j| \geq  2$. }
\end{equation}

(i)  If $A=B$, then $s\geq |A|+|B|+f-1$. 

(ii) If $B\neq A$, then $s \geq |A| + |B| + f$ with equality only if $A\subset B$ or $B\subset A$.
\end{lemma}

\begin{proof} Let $Q'=v'_1, e'_1,v'_2, \ldots, e'_{s'-1},v'_{s'}$ be a path obtained from $Q$
by iteratively contracting $f$ edges of $F$. In particular, $s'=s-f$.
Since $A$ and $B$ are both vertex-disjoint from $F$, each $v_i \in A \cup B$ was unaffected by the edge contractions and hence still exists as some $v'_{i'}$ in $Q'$. Moreover,~\eqref{dis2} still holds for $A$ and $B$ in $Q'$.

So,  Lemma~\ref{ver-new} for $q=2$ applied to $A,B$ and $Q'$ yields that if $A=B$, then $s'\geq |A|+|B|-1$, and 
if $B\neq A$, then $s' \geq |A| + |B| $ with equality only if $A\subset B$ or $B\subset A$. Since $s'=s-f$, this proves our lemma.
\end{proof}

\section{Setup for  Theorem~\ref{main}}

Bounds of the theorem differ for $r\leq n/2$ and $r>n/2$. Naturally, the proofs also will be different, but they will have similar structure. In both proofs, for given vertices $x,y$ in an $r$-graph $H$ we attempt to find a hamiltonian
$x,y$-path.
Both proofs will have three steps.

In Step 1 we construct an $x,y$-path $Q$ with at least $\max \{\lceil \frac{n+2}{2}\rceil,r+1\}$ vertices.

Then we consider  pairs $(Q, P)$ of  vertex-disjoint paths in $H$ such that $Q$ is an $x, y$-path.
We will say that such a pair $(Q,P)$ {\em is better} than a similar pair $(Q',P')$ if 
\begin{enumerate}
    \item[(i)] $|E(Q)|>|E(Q')|$, or
    \item[(ii)] $|E(Q)|=|E(Q')|$ and $|E(P)|>|E(P')|$,
or
    \item[(iii)] $|E(Q)|=|E(Q')|$, $|E(P)|=|E(P')|$ and the total number of vertices in $V(P)$ in the edges in $Q$ (counted with multiplicities) is greater than the total number of vertices in $V(P')$ in the edges in $Q'$.
\end{enumerate}

We consider best pairs and study their properties. Some properties will be proven in the next subsection. Using these properties together with the lemmas on graph paths from the previous section  in Step 2 we show that the path $P$ in a best pair cannot have exactly one vertex.
In the final Step 3 we handle all cases when $P$ has at least two vertices.

\medskip
Below we  assume that
 $(Q,P)$ is a best pair, $Q = v_1, e_1, \ldots, e_{s-1}, v_s$,  and $P = u_1, f_1, \ldots, f_{\ell-1}, u_{\ell}$. 

We consider three subhypergraphs, $H_Q,H_P$ and $H'$ of $H$ with the same vertex set $V(H)$:  $E(H_Q)=\{e_1,\ldots,e_{s-1}\}$, $E(H_P)=\{f_1,\ldots,f_{\ell-1}\}$ and  $E(H')=E(H)-E(H_P)-E(H_Q)$. 
By definition, the edge sets of these three subhypergraphs form a partition of the edge set of $H$.
For a hypergraph $F$ and a vertex $u$, we denote by $N_{F}(u) = \{v \in V(F): \{u,v\} \subseteq e \text{ for some } e \in F\}$. 
For $i \in \{1, \ell\}$, set $B_i=\{e_j \in E(Q): u_i \in e_j\}$ and $b_i = |B_i|$. 

\subsection{Claims on best pairs}
The  claims below apply to all best pairs $(Q,P)$, regardless of the uniformity $r$. 

\begin{claim}\label{noconsecutive}
In a best pair $(Q,P)$, $N_{H'}(u_1)$ cannot contain a pair of vertices that are consecutive in $Q$.  
\end{claim}

\begin{proof}
Suppose toward a contradiction that $v_i, v_{i+1}$ are contained in edges of $H'$ with $u_1$. Let $e, e' \in E(H')$ be such that $u_1, v_i \in e$ and $u_1, v_{i+1} \in e'$. If $e \neq e'$, then replacing $e_i$ with $e,u_1, e'$ gives a longer $x,y$-path than $Q$, a contradiction. Thus we may assume $e = e'$.

If there is $1 \leq j \leq \ell$ such that $u_j \in e_i$, then by replacing the path $v_i, e_i, v_{i+1}$ in $Q$ with the longer path $v_i, e, u_1, f_1, u_2, \ldots, f_{j-1}, u_j, e_i, v_{i+1}$, we obtain a longer $x,y$-path than $Q$. Thus $e_i \cap V(P) = \emptyset$. Then replacing $e_i$ with $e$ in $Q$ gives a path $Q'$ with $(Q',P)$ better than $(Q,P)$ by criterion (iii).
\end{proof}

Symmetrically, the claim holds for $u_\ell$ in place of $u_1$.

\begin{claim}\label{notneighbor}
For any $u \notin V(Q)$, if $u \in e_i$, then $v_i, v_{i+1} \notin N_{H - H_Q}(u)$. 
\end{claim}

\begin{proof}
Suppose $v_i \in N_{H - H_Q}(u)$, and let $e \in E(H)-E(H_Q)$ be such that $\{u, v_i\} \subseteq e$. Then we can find a longer cycle by replacing $e_i$ with $e, u, e_i$, a contradiction to our choice of $Q$. A similar argument holds for $v_{i+1}$. 
\end{proof}

\begin{claim}\label{distance}For every $e_i \in B_1, e_j \in B_\ell$ either $i=j$ or $|i - j| \geq \ell$.\end{claim}

\begin{proof}
Suppose there exists $e_i \in B_1, e_j \in B_\ell$ such that without loss of generality $j > i$ and $j - i \leq \ell-1$. Then the  path obtained by replacing $v_i, e_i, \ldots, e_{j}, v_{j+1}$ in $Q$ with $v_i, e_i, u_1, f_1, \ldots, f_{\ell-1}, u_\ell, e_j, v_{j+1}$ has  $|V(Q)| - (i-j) + \ell > |V(Q)|$ vertices, a contradiction.
\end{proof}

\begin{claim}\label{distance2}If there exists distinct edges $e, f \in E(H')$ such that $\{u_1, v_i\} \subset e$ and $\{u_\ell, v_j\} \subset f$, then  $|i - j| \geq \ell+1$. \end{claim}

\begin{proof}If $|i-j| \leq \ell$, replace the subpath in $Q$ from $v_i$ to $v_j$ with the path $v_i, e, P, f, v_j$ to get a longer $x,y$-path.\end{proof}

\begin{claim}\label{distance3}For every $v_i \in N_{H'}(u_1)$ and $e_j \in B_\ell$, 
if $i \leq j$ then $j-i \geq \ell$ and if $i > j$ then $i-j \geq \ell+1$.
\end{claim}

\begin{proof}Let $e \in E(H')$ contain $v_i$ and $u_1$. If $i\leq j$, let $Q'$ be the path obtained by replacing the segment $v_i, e_i, \ldots, e_j, v_{j+1}$ in $Q$ with the path $v_i, e, u_1, P, u_\ell, e_j, v_{j+1}$. If $i > j$, let $Q'$ be obtained from $Q$ by replacing $v_j, e_j, \ldots, v_i$ with $v_j, e_j, u_\ell, f_{\ell-1}, \ldots, f_1, u_1, e, v_i$. In the first case, $|V(Q')| = |V(Q)| - (j-i) + \ell$, and in the second case $|V(Q')| = |V(Q)| - (i-(j+1)) + \ell$. The Claim follows since $|V(Q)| \geq |V(Q')|$ by the choice of $(Q,P)$. 
\end{proof}

\begin{claim}\label{forbidden}
For any $e \in E(H')$, if $v_i, v_j \in e$, then at most one of $e_i, e_j$ is in $B_1$ and at most one of $e_{i-1}, e_{j-1}$ is in $B_1$. 
\end{claim}

\begin{proof} If $e_{i-1}, e_{j-1} \in B_1$, then we get a longer path \[
 v_1, e_1, v_2, \dots, v_{i-1}, e_{i-1}, u_1, e_{j-1}, v_{j-1}, e_{j-2}, v_{j-2}, \dots, v_i, e, v_j, e_j, v_{j+1}, \dots, v_s.
 \]
 
 The argument for $e_i, e_j$ is similar.
\end{proof}

\begin{claim}\label{forbidden2}
Let $B_1^-=\{v_i: e_i\in B_1\}$ and $B_1^+=\{v_{i+1}: e_i\in B_1\}$.

 (i) For any edge $e \in E(H')$, $b_1 \leq s -|e \cap V(Q)|+1$ with equality only if $B_1 = \{e_i, e_{i+1}, \ldots, e_{j}\}$ for some $i<j$ and $e \cap V(Q)= \{v_1, \ldots, v_i\} \cup \{v_{j+1}, \ldots, v_s\}$. 

(ii) $b_1 \leq s-1 - |N_{H'}(u_1)\cap V(Q)|$ with equality only if $B_1$ is a set of $b_1$ consecutive edges in $Q$ and $N_{H'}(u_1) \cap V(Q) = V(Q) - (B_1^- \cup B_1^+)$.\end{claim}

\begin{proof}
 Set $e' = e \cap V(Q)$. By Claim~\ref{forbidden}, $|B_1^- \cap e| \leq 1$.
  Hence
  $$b_1 -1 + |e \cap V(Q)|-1 \leq |B_1|-1+|e' - \{v_s\}|\leq |\{v_1, \ldots, v_{s-1}\}|=s-1,$$ i.e., $b_1 \leq s - |e'| + 1$
  with equality only if $v_s\in e$, $|B_1^- \cap e|=1$, and $e\cup B_1^-=V(Q)$. Symmetrically, $v_1\in e$,
  $|e \cap B_1^+|=1$, and $V(Q)\setminus e\subset B_1^+$. So,
   $V(Q)\setminus e\subseteq B_1^-\cap B_1^+$, and therefore $|B_1^-\cap B_1^+|\geq s-|e'| = b_1 -1$.
   This means, the symmetric difference of $B_1^-$ and $B_1^+$ has only two vertices.
    For this, the $b_1$ edges in $B_1$ must be consecutive on $Q$, and $e' = V(Q) - (B_1^- \cap B_2^+)$. This proves (i).
    
For (ii), by Claim~\ref{notneighbor} $(B_1^- \cup B_1^+) \cap  (N_{H'}(u_1)\cap V(Q)) = \emptyset$. Therefore $|B_1^- \cup B_1^+| \leq |V(Q)| -|N_{H'}(u_1)\cap V(Q)| = s - |N_{H'}(u_1)\cap V(Q)|$. We have $|B_1^- \cup B_1^+| \geq b_1 + 1$ with equality only if $B_1$ is a set of consecutive edges in $Q$. These inequalities together give our result.
\end{proof}

\begin{cor}\label{b1bound}
When $\ell = 1$, $b_1 \leq n-3$. If in addition,  $|N_{H'}(u_1) \cap V(Q)| \geq 2$, then $b_1 \leq n-4.$
\end{cor}

\begin{proof}
%

Since $|E(Q)| < n-1$, there must be at least one edge $e \in E(H')$, and since $\ell = 1$, $e$ contains at least $r-1 \geq 2$ vertices in $Q$ if $u_1 \in e$, and at least $r \geq 3$ otherwise. By Claim~\ref{forbidden2}, $b_1 \leq s-2\leq n-3$. 

The second part follows from Claim~\ref{forbidden2} (ii).
\end{proof}

\section{Finding a longish $x,y$-path}

In this section, we will show that there exists an $x,y$-path of length at least $\max\{n/2 + 1, r+1\}$. 

\begin{lemma}\label{longpath}Let $3 \leq r \leq n/2$, and let $H$ be an $n$-vertex $r$-graph. Let $x,y\in V(H)$. If $\delta(H) \geq {\lfloor n/2 \rfloor \choose r-1}$, then $H$ contains an $x,y$-path with at least $n/2 +1$ vertices.\end{lemma}
\begin{proof}The restrictions on $\delta(H)$
 in Theorem~\ref{oldmain} are not stronger than  in Theorem~\ref{main}.
 So, by Theorem~\ref{oldmain}, $H$ contains a hamiltonian cycle $C$ in $H$. The longer of the two $x,y$-paths along $C$ has  at least $n/2+1$ vertices. 
\end{proof}
 
  For $r>n/2$, we need much more effort, see below.

\begin{lemma}\label{r+1path} Let $n \geq r >n/2$, and let $H$ be an $n$-vertex $r$-graph. Let $x,y\in V(H)$. If $\delta(H) \geq r-1$, then $H$ contains an $x,y$-path with at least $r+1$ vertices.\end{lemma}

 \begin{proof}
 We will first show that there exists some $x,y$-path in $H$. If there exists an edge $e \in E(H)$ with $\{x,y\} \subseteq e$, then we are done. Otherwise since $r > n/2$, any two edges $e,f \in E(H)$ such that $ x\in e, y \in f$ have a common vertex, say $v \in e \cap f$. Then $x,e,v,f,y$ is an $x,y$-path  in $H$.
 
 Now let $Q = v_1, e_1, \ldots, e_{s-1}, v_s$ be a longest $x,y$-path in $H$ (so $x = v_1, y = v_s$). Moreover, choose $Q$ so that if $\{v_1, \ldots, v_s\} \in E(H)$, then this edge is used in $Q$. Suppose  $s\leq r$. 
 
  Construct a new hypergraph $\hat{H}$ as follows: $V(\hat{H}) = V(H) - V(Q)$, and $E(\hat{H}) = \{e \cap V(\hat{H}): e \in E(H) - E(Q)\}$. Note that $\hat{H}$ is not necessarily a uniform hypergraph. We have a mapping from the edges of $H$ to the edges of $\hat{H}$ given by $e \mapsto e - V(Q)$ (which is not necessarily one-to-one). 

Let $D_1, D_2, \ldots, D_q$ be the vertex sets of the connected components of $\hat{H}$. For $1\leq j\leq q$, let $d_j = |\{e_i \in E(Q): e_i \cap D_j\neq \emptyset\}|$.  
Since $|V(Q)| \leq r$, at most one edge $e_i \in E(Q)$ may be contained in $V(Q)$. It follows that 

\begin{equation}\label{djsum}\sum_{i=1}^q d_i \geq |E(Q)|-1 = s-2.\end{equation}

\begin{claim}\label{diffcomp} For any $1\leq j\leq q$, if $e_i\cap D_j\neq \emptyset$, then the edges of $E(H) - E(Q)$ containing  $v_i$ or $v_{i+1}$ cannot intersect $D_j$.
\end{claim}

\begin{proof}Let $v \in D_j \cap e_i$. Suppose $v_i \in h \in E(H)-E(Q)$ and $u\in h \cap D_j $. 
  Then  $\hat{H}$ contains a $u,v$-path  which we can lift to a $u,v$-path $P$ in $H$ that avoids $E(Q)$. If $ h \notin E(P)$, then by replacing the segment $v_i, e_i, v_{i+1}$ in $Q$ with the path $v_i, h,  P, e_i, v_{i+1}$, we obtain a longer $x,y$-path. Otherwise let $P'$ be the subpath of $P$ starting with $h$. Then  we replace $v_i, e_i, v_{i+1}$ with $v_i, P', e_i, v_{i+1}$ to get a longer path. The argument  for $v_{i+1}$ is similar.
\end{proof}

%

\begin{claim}\label{consec}For any $1\leq j\leq q$ and  any $1\leq i \leq s-1$, there are no distinct edges $e, f \in E(H) - E(Q)$ such that $e$ and $f$ intersect $D_j$, $v_i \in e$, and $v_{i+1} \in f$. \end{claim}

\begin{proof}
Let $P'$ be a shortest path in $\hat{H}$ from  $e\cap D_j$ to $f\cap D_j$. Lift $P'$ to a path $P$ in $H$ which avoids $E(Q)$. By the minimality of $P'$, $e \notin E(P)$ and $f \notin E(P)$.
Then we may replace the segment $v_i, e_i, v_{i+1}$ in $Q$ with $v_i, e, P, f, v_{i+1}$ to get a longer $x,y$-path.
\end{proof}

\begin{claim}\label{deg2}For any $1\leq j\leq q$, if at least 2 edges in $E(H) - E(Q)$ intersect $D_j$, then \[|D_j| \geq r - \lceil (s-d_j)/2 \rceil + 1.\]\end{claim}
\begin{proof}
Suppose $|D_j| \leq r- \lceil (s-d_j)/2 \rceil$, and let $e, g \in E(H) - E(Q)$ be distinct edges that intersect $D_j$. 
 Let $A=e \cap V(Q)$, $B=g \cap V(Q)$, and $F=\{\{v_i,v_{i+1}\}:\, D_j\cap e_i \neq \emptyset\}$.
 By definition, $|F|=d_j$, and
each of  $A$ and $B$ has at least $r - |D_j| \geq \lceil (s-d_j)/2 \rceil$ vertices.

  By Claim~\ref{diffcomp}, $A$ and $B$ are disjoint from all pairs in $F$. By Claim~\ref{consec},~\eqref{dis2} holds.
So
    Lemma~\ref{consecpath2} together with the lower bounds on $|A|$ and $|B|$ imply that if $A\neq B$, then 
    \begin{equation}\label{shortp}
    s\geq |A|+|B|+d_j\geq 2\frac{s-d_j}{2}+d_j=s,
    \end{equation}
     with equality only if
   $A\subset B$ or $B\subset A$. But if $A\subset B$ or $B\subset A$ and $A\neq B$, then $|A|+|B|\geq 1+2\frac{s-d_j}{2}$. Hence, if $A\neq B$, then in the RHS of~\eqref{shortp} we get at least $s+1$, a contradiction.
  
  Thus $A=B$. Since $e$ and $g$ are distinct but coincide on $Q$, $e \cap D_j$ and $g \cap D_j$ are distinct sets each with at least $r-\lceil (s-d_j)/2 \rceil$ vertices. It follows that $|D_j| \geq r  - \lceil (s-d_j)/2 \rceil + 1$.
\end{proof}

\begin{claim}\label{deg1}For any $1\leq j\leq q$, if exactly one edge in $E(H) - E(Q)$ intersects $D_j$, then $|D_j| \geq r$. \end{claim}
\begin{proof}
Suppose $|D_j| \leq r - 1$ and $e$ is the unique edge in $E(H) - E(Q)$ that intersects $D_j$. Then by the definition of $\hat{H}$, $D_j=e-V(Q)$. Let $v\in D_j$.
 Since $|e| = r$, $e$  contains at least one vertex $v_i$ in $Q$. By symmetry, we may suppose $i < s$. In order to have $d(v) \geq r-1$, $v$ must belong to at least $r-2$ edges of $E(Q)$. By Claim~\ref{diffcomp},
 none of these at least $r-2$ edges is $e_{i-1}$ or $e_i$. This is possible only if $s=r$, $e\cap V(Q)=\{v_1\}$
 and $v\in e_2\cap e_3\cap\ldots \cap e_{s-1}$.
 This implies $|D_j| = r-1$ and each vertex in $D_j$ belongs to $e_2$ by symmetry. But then $\{v_2, v_3\} \cup D_j \subseteq e_2$,
 contradicting the fact that $|e_2| = r$.
\end{proof}

\begin{claim}\label{deg0}For any $1\leq j\leq q$, at least one edge in $E(H) - E(Q)$ intersects $D_j$.\end{claim}

\begin{proof}
Suppose not. By the definition of $\hat{H}$, $D_j$ is a single vertex, say $v$. Since $d(v) \geq r-1$, $v$ must belong to at least $r-1$ edges of $Q$, which is only possible if $|V(Q)| = r$. In this case $v$ is contained in all edges of $Q$.  By the choice of $Q$, we have  $\{v_1, \ldots, v_s\} \notin E(H)$.

Since $|E(H)| \geq n-1 > r-1$, there exists an edge  $g\in E(H) - E(Q)$. By the choice of $Q$, $g$ intersects some $D_h$. If $|D_h| \geq r-1$, then $|V(H)| \geq |V(Q)|+ |D_h| + |D_j| \geq r + r-1 + 1 > n$, a contradiction. In particular, by Claim~\ref{deg1}, this implies that at least two edges in $E(H) - E(Q)$ intersect $D_h$. We claim that for each such edge $e$, 
\begin{equation}\label{v1vs}e \cap V(Q) \subseteq \{v_1, v_s\}.\end{equation} 
Suppose this is not the case. Then since $|D_h| \leq r-2$, there exists a pair $\{v_i, v_{i'}\} \neq \{v_1, v_s\}$ and edges $e, f \in E(H) - E(Q)$ such that $e$ and $f$ intersect $D_h$, $v_i \in e$, and $v_{i'} \in f$. Without loss of generality, we may assume $i < i' < s$. Let $P$ be a $v_i, v_{i'}$-path in $H$  avoiding $E(Q)$ (it could be the case that $P$ contains only one edge). Then 
\[v_1, \ldots, v_i,P, v_{i'}, e_{i'-1}, \ldots, v_{i+1}, e_i, v, e_{i'}, v_{i'+1}, \ldots, e_{s-1}, v_s\] is a longer $x,y$-path in $H$.  Therefore~\eqref{v1vs} holds. Since at least two edges in $E(H) - E(Q)$ intersect $D_h$, $|D_h| \geq  r-1$, a contradiction.
\end{proof}

\begin{claim}\label{2comps}$\hat{H}$ has at least 2 components.
\end{claim}
\begin{proof}
Suppose $q=1$.  
 If $|V(Q)| \leq r-1$, then each edge $e_i$ intersects $D_1$ and each $v_i \in V(Q)$ is contained in an edge $h \in E(H) - E(Q)$, and $h$ must also intersect $D_1$ since $|h| > |V(Q)|$, contradicting Claim~\ref{diffcomp}. So we may assume $|V(Q)| = r$.

By Claim~\ref{deg0}, some edge $h \in E(H) - E(Q)$  intersects $D_1$. If $h \subset D_1$, then $|V(H)| \geq |V(D_1)| + |V(Q)| \geq r + r > n$, a contradiction. For each $v_i \in h \cap V(Q)$, each of $e_i$ and $e_{i-1}$ must be contained in $V(Q)$. As $r = V(Q)$, only one such edge in $Q$ can satisfy this. Hence without loss of generality, we may assume $h\cap V(Q) \subseteq \{v_1\}$ and $e_1 = V(Q)$. It follows that $|D_1| \geq r-1$. If $|D_1| \geq r$, then again we get $|V(H)| > n$.

Hence, the last possibility is that
 $|D_1| = r-1$ and $h \cap V(Q) = \{v_1\}$. In particular, by Claim~\ref{deg1}, some other edge $h' \in E(H) - E(Q)$  intersects $D_1$. Since $s=r\geq 3$, $e_1\neq e_{s-1}$. So by the same argument as for $h$, we have $h'\cap V(Q) = \{v_1\}$.
 Since $h'\neq h$ and $D_1\supseteq h\cup h'-\{v_1\}$, we get $|D_1|>r-1$, a contradiction.
\end{proof}

Now we are ready to finish the proof of the lemma. By Claims~\ref{deg0},~\ref{deg1} and~\ref{deg2}, $|V(D_j)| \geq r - \lceil (s-d_j)/2 \rceil + 1$ for all $j$. Therefore

\[|V(H)| \geq  |V(Q)| + \sum_{j=1}^q (r - \lceil (s-d_j)/2 \rceil + 1)
\geq  s + q(r - \frac{s+1}{2} + 1) + \sum_{i=1}^q \frac{d_j}{2}.
\]

Since 
 $r\geq s$, the quantity $q(r - \frac{s+1}{2} + 1)$ is minimized when $q = 2$. By~\eqref{djsum}, 
\begin{eqnarray*}
|V(H)| &\geq &  s + 2(r - \frac{s+1}{2} + 1) + \sum_{i=1}^q \frac{d_j}{2}\\
&\geq &  s + 2r - (s+1) + 2 + (s-2)/2\\ 
& = & 2r + s/2\\
&>& n,
\end{eqnarray*}
a contradiction.
 \end{proof}

\section{Proof of Theorem~\ref{main} for $r \leq n/2$}

In the next two sections, we set $t' = \lfloor n/2 \rfloor$ and consider a best pair $(Q,P)$ with $Q = v_1, e_1, \dots, e_{s-1}, v_s$ and $P = u_1, f_1, \dots, f_{\ell-1}, u_\ell$. By Lemma~\ref{longpath}, $s \geq t'+1$ if $s$ is even and $s \geq t'+2$ if $s$ is odd. In both cases we get $\ell \leq n-s \leq t'-1$ and $s \geq n/2+1$. Recall that for $i \in \{1,\ell\}$, $B_i = \{e_j: u_i \in e_j\}$, and $b_i = |B_i|$.

\subsection{Finding a nontrivial path $P$}\label{step2}

\begin{lemma}\label{ell1} In a best pair $(Q,P)$, $|V(P)| \geq 2$.\end{lemma}

\begin{proof}
Suppose that $|V(P)|=\ell=1$, i.e., $P = u_1$. Then $s \leq n-1$. By condition (ii) of $(Q,P)$ being a best pair,  every edge of $H'$ contains at most one vertex outside $Q$. 

 Claims~\ref{noconsecutive},~\ref{notneighbor} and Lemma~\ref{indep2} imply that 
 $|N_{H'}(u_1)| \leq \lceil (s-b_1)/2 \rceil$. Therefore

\begin{equation}\label{degbound}
1 + {t' \choose r-1} \leq d_H(u_1) \leq b_1 + {\lceil (s-b_1)/2 \rceil \choose r-1} \leq b_1 + {\lceil (n-1-b_1)/2 \rceil \choose r-1}.
\end{equation}

\medskip
{\bf Case 1}: $b_1 = 0$. By~\eqref{degbound}, $1 + {t' \choose r-1}  \leq {\lceil (n-1)/2 \rceil \choose r-1} = {t' \choose r-1}$, a contradiction. 

\medskip
{\bf Case 2}: $b_1 = 1$. Again by~\eqref{degbound},  $1 + {t' \choose r-1}  \leq 1 + {\lceil (s-1)/2 \rceil \choose r-1} \leq 1+ {\lceil (n-2)/2 \rceil \choose r-1}$. If $n$ is even, we immediately obtain a contradiction. If $n$ is odd, then we reach a contradiction when $s < n-1$. So suppose $n$ is odd, $s = n-1$, $|N_{H'}(u_1)| = \lceil (s-b_1)/2 \rceil = s/2=t'$, and $u_1$ is contained in all ${t' \choose r-1}$ possible edges within $N_{H'}(u_1) \cup \{u_1\}$. 

Consider the unique edge $e_i$ of $Q$ containing $u_1$. Then $|N_{H'}(u_1)| \leq \lceil (i-1)/2 \rceil + \lceil (n-1-(i+1))/2 \rceil$ by Claim~\ref{noconsecutive} and Claim~\ref{notneighbor}. If $i$ is odd, then this gives  $|N_{H'}(u_1)| \leq (i-1)/2 + (n-i-2)/2 = (n-3)/2$, a contradiction. Thus, $i$ is even and $X := N_{H'}(u_1) = \{v_1, v_3, \dots, v_{i-1}, v_{i+2}, v_{i+4}, \dots, v_s\}$.

Replacing $e_{i-1}$ in $Q$ with the edge $e \in E(H')$ containing $u_1, v_{i-1}$ and replacing $v_i$ with $u_1$ creates a new path $Q'$ which only misses $v_i$. Since $(Q,P)$ is a best pair, by condition (iii) of choosing a best pair, $e_i$ and $e_{i-1}$ can be the only edges of $Q$ which contain $v_i$ and in fact $(Q',v_i)$ is also a best pair. Thus applying the same arguments to $v_i$ and $Q'$ as we did to $u_1$ and $Q$, we obtain that $N_{H-Q'}(v_i) = X$. Notice that we can apply a symmetric argument to $v_{i+1}$ and corresponding path $Q''$ to get $N_{H-Q''}(v_{i+1}) = X$.

We will find an edge $g \neq e_i$ with $|g-X| \geq 2$ and $|g\cap \{v_2, v_4, \dots, v_{i-2}, v_{i+3}, v_{i+5}, \dots, v_{s-1}\}| \geq 1$, and then we will use $g$ to find the desired hamiltonian path. Choose $v_j \notin e_i$ with $v_j \notin X$, which exists because $|(X \cup e_i) \cap V(Q)| \leq |X| + r-1 \leq 2t'-1 = s-1$. Since $d_H(v_j) > {t' \choose r-1}$ and $|X| = t'$, there is an edge $g$ containing $v_j$ and some vertex outside $X$. Since $v_j \notin X$ and $v_j \notin e_i$, that vertex cannot be $u_1$ and must instead be some $v_k \in V(H) - (X \cup \{u_1\}) =  V(Q)-X$. Suppose without loss of generality that $j <k$.

{\bf Case 2.1}: $g\in E(H')$. Since $v_j$ is in neither $X$ nor $e_i$, $v_{j-1} \in X$. Thus let $f \in E(H')$ be such that $v_{j-1}, u_1 \in f$. Similarly, since $v_k \notin X$, we have $v_{k-1} \in X$ unless $k = i+1$, which we handle separately. Let $f' \in E(H')$ be such that $v_{k-1}, u_1 \in f'$, and observe that we can choose edges such that $f, f'$ are distinct because $u_1$ is in ${t'-1 \choose r-2} \geq 2$ edges with each vertex in $X$. Thus if $j< k$, we have the hamiltonian path 

\[
v_1, e_1, v_2, \dots, v_{j-1}, f, u_1, f', v_{k-1}, e_{k-2}, v_{k-2}, \dots, v_j, g, v_k, e_k, v_{k+1}, \dots, v_s.
\]

A similar path can be found for $j > k$ by symmetry. In the case $k=i+1$, replace $f'$ in the above path with $e_i$ to obtain the desired hamiltonian path. 

{\bf Case 2.2}: $g = e_m \in E(Q)$. Since $g \neq e_i$, we may assume by symmetry that $v_m \in X$, unless $m=k$. Let $f$ be as in the previous case, and let $f' \in E(H')$ be such that $v_m, u_1 \in f'$. 

Thus for $j <m$ we have the hamiltonian path

\[
v_1, e_1, v_2, \dots, v_{j-1}, f, u_1, f', v_{m}, e_{m-1}, v_{m-1}, \dots, v_j, g, v_{m+1}, e_{m+1}, v_{m+2}, \dots, v_s
\]
(and similar for $j >m$). If $m=k$, then $v_{k+1} \in X$, so we let $f'', f''' \in E(H')$ be such that $v_{j+1}, u_1 \in f''$ and $v_{k+1}, u_1 \in f'''$. Then 

\[
v_1, e_1, v_2, \dots, v_{j}, g, v_k, e_{k-1}, v_{k-1}, \dots, v_{j+1}, f'', u_1, f''', v_{k+1}, e_{k+1}, v_{k+2} , \dots, v_s
\]

is hamilitonian if $j < k$ and we can find a similar path for $j > k$, ending the proof of Case 2.


\medskip
{\bf Case 3}: $b_1 \geq 2$. Then by (\ref{degbound})
\[
1 + {t' \choose r-1} \leq  b_1 + {\lceil (n-1-b_1)/2 \rceil \choose r-1} \leq b_1 + {\lceil (n-3)/2 \rceil \choose r-1} = b_1 + {t' -1 \choose r-1}.
\]

Hence $1 + {t' \choose r-1} - {t' - 1 \choose r-1} = 1 + {t'-1 \choose r-2} \leq b_1 \leq n-3$ by Corollary~\ref{b1bound}. If $2 \leq r-2 \leq t'-3$, then we have $n-4 \geq {t'-1 \choose r-2} \geq {t'-1 \choose 2}$, a contradiction when $n \geq 12$. For $n \leq 11$, it is straightforward to check that $1 + {t' \choose r-1} >  b_1 + {\lceil (n-1-b_1)/2 \rceil \choose r-1}$ in all cases.

For $r=3$, we have $t' = 1 + {t'-1 \choose r-2} \leq b_1$, so 

\[
1 + {t' \choose 2} \leq b_1 + {\lceil (n-1-b_1)/2 \rceil \choose 2} \leq b_1 + {t' -\lfloor b_1/2 \rfloor \choose 2} \leq n-3 + {\lceil t'/2 \rceil \choose 2}.
\]

This gives a contradiction when $n \geq 12$. For $n \leq 11$, it is straightforward to check that $1 + {t' \choose 2} >  b_1 + {\lceil (n-1-b_1)/2 \rceil \choose 2}$ in all cases except $n = 7$, $b_1 \in \{3,4\}$, which will be considered with the case $r=t'=3$. 

For $r = t'$, by (\ref{degbound}) we have
\[
1+t' = 1 + {t' \choose r-1} \leq b_1 + {\lceil (s-b_1)/2 \rceil \choose r-1} \leq b_1 + 1,
\]

since $s \leq n-1$ and $b_1 \geq 2$. Thus $b_1 \geq t'$. For $n \geq 8$, we have 
\[
|N_{H'}(u_1)| \leq \lceil (s-b_1)/2 \rceil \leq \lceil (n-1-t')/2 \rceil \leq \lceil (n-5)/2 \rceil = t'-2 < r-1.
\] This also holds for $n < 8$ if $s \leq n-2$, so we will handle the case $n <8, s=n-1$ separately at the end of this subsection.

Since each edge in $H'$ containing $u_1$  contains $r-1$ other vertices, $|N_{H'}(u_1)| < r-1$ gives that $|N_{H'}(u_1)|=0$ and hence $1 + t' \leq b_1$. Notice also that $n\leq \delta(H) \frac{n}{t'}  \leq |E(H)| = |E(Q)| + |E(H')| \leq n-2 +|E(H')|$, so there are at least 2 edges in $E(H')$.

\medskip
{\bf Case 3.1}: There exists $e \in E(H')$ with $e \subseteq V(Q)$. By Claim~\ref{forbidden2}, $b_1 \leq s -t'+1 \leq n-1 -t' + 1 \leq t'+1$ with equality only if there exists $i<j$ such that $B_1 = \{e_i, \ldots, e_j\}$ and $e = \{v_1, \ldots, v_i\} \cup \{v_{j+1}, \ldots, e_s\}$. Without loss of generality, $i \geq 2$. Then we can replace $e_1$ in $Q$ with $e$ to obtain another $x,y$-path $Q'$ such that $(Q', P)$ is also a best pair. As $B_1$ does not change for this new pair, $e_1$ must play the old role of $e$, i.e., $e_1 = \{v_1, \ldots, v_{i}\} \cup \{v_{j+1}, v_{j+3}, \ldots, v_s\}$, but then $e = e_1$, a contradiction.


\medskip
{\bf Case 3.2}: Every edge in $H'$ contains exactly one vertex outside of $V(Q)$. Since $E(H') \neq \emptyset$ and $u_1$ is contained only in edges of $Q$, there must be at least one additional vertex outside of $Q$ and hence $s \leq n-2$. Let $e \in E(H')$. Because Case 3.1 does not hold, $|e - V(Q)|\leq n-1 -s$, so $|e \cap V(Q)| \geq t' - (n-1-s)$ with equality only if $e \cup V(Q) \cup \{u_1\} = V(H)$. By Claim~\ref{forbidden2}, $1 + t' \leq b_1 \leq s - |e \cap V(Q)| + 1 \leq s - (t'-(n-1-s))+1\leq +n-t'$, so $2t' +1 \leq n$. We get a contradiction unless the ``equality" part of Claim~\ref{forbidden2}(i) holds. Then as in the previous subcase, there exists $i<j$ such that $e \cap V(Q)= \{v_1, \ldots, v_i\} \cup \{v_j, \ldots, v_s\}$. Moreover, since the choice of $e \in E(H')$ was arbitary, for each $e' \neq e$ in $E(H')$, $e' \cap V(Q) = e \cap V(Q)$. But since $V(H) = e \cup V(Q) \cup \{u_1\}$ and $u_1 \notin e'$, $e' - V(Q) = e-V(Q)$, hence $e = e'$, a contradiction.

\medskip

Finally we handle the cases $6 \leq n \leq 7$, $r=t'=3$, $s= n-1$, and $b_1 \in \{3,4\}$. The average degree of $H$ is \[\sum_{v \in V(H)} d(v)/n = 3|E(H)|/n \geq \delta(H) \geq 4,\] so $|E(H)| \geq \lceil 4n/3 \rceil$ which is equal to $8$ when $n=6$ and $10$ when $n=7$. In either case, there exists at least 3 edges in $H'$.

 We will first show that $B_1$ is a set of $b_1$ consecutive edges in $Q$. If $u_1$ is not contained in any edges in $H'$, then $b_1 \geq \delta(H) \geq 4$. Otherwise if $u_1$ belongs to an edge $h$ of $H'$, then Claim~\ref{forbidden2}(ii) implies $b_1 = 3$ and $n=7$. In both cases, the ``equality" part of Claim~\ref{forbidden2} implies $B_1 = \{e_i, e_{i+1}, \ldots, e_{i+b_1-1}\}$ for some $i$. 
 
If $u_1$ is not contained in any edges in $H'$, then for any $e \in E(H'), e  = \{v_1, \ldots, v_i\} \cup \{v_{i+b_1}, \ldots, v_s\}$. But this holds for all edges in $H'$, a contradiction.
Now suppose $n=7$, $b_1 = 3$ and $u_1$ is contained in an edge $e$ of $H'$. Since $N_{H'}(u_1)$ contains no consecutive vertices and is disjoint from $\{v_i,v_{i+1},v_{i+2},v_{i+3}\}$, we  have $i=2$, and $e = \{v_1, v_6, u_1\}$. In particular, $d_{H'}(u_1)=1$.

Let $E''$ be the set of edges in $H$ not containing $u_1$. Since $|E''|=|E(H)|-d_H(u_1)\geq 10-4=6$, some edge $g\in E''$ does not contain $\{v_1,v_6\}$. By symmetry, we may assume $v_6\notin g$. If $g=e_1=\{v_1,v_2,v_h\}$, then we have a longer $v_1,v_6$-path $v_1,e,u_1,e_{h-1},v_{h-1},e_{h-2},\ldots,v_2,e_1,v_h,e_h,v_{h+1},\ldots,v_6$.

Otherwise, $g\in E(H')$. So, by Claim~\ref{forbidden}, $|g\cap \{v_2,v_3,v_4\}|\leq 1$ and $|g\cap \{v_3,v_4,v_5\}|\leq 1$.
This is possible only if $g= \{v_1,v_2,v_5\}$. Then we have  $v_1,v_6$-path $v_1,e,u_1,e_{4},v_{4},e_{3},v_3,e_2,v_2,g,v_5,e_5,v_6$, a contradiction.
%
%
%
%
\end{proof}

\subsection{Finishing the proof of Theorem~\ref{main} for $r\leq n/2$}

\begin{proof}[Proof of Theorem~\ref{main} for $r \leq n/2$.]
Consider a best pair $(Q,P)$ with $Q = v_1, e_1, \dots, e_{s-1}, v_s$ and $P = u_1, f_1, \dots, f_{\ell-1}, u_\ell$.

By symmetry, we may assume $b_\ell = |B_\ell| \geq |B_1| = b_1$. By Lemma~\ref{ell1}, $\ell \geq 2$. Recall also that $s \geq n/2 + 1 \geq t'+1$ and $\ell \leq t'-1$. 


 By Claim~\ref{distance} and Lemma~\ref{ed-new}, either 
\begin{equation}\label{bineq1}
    s\geq  b_1 +b_\ell+\ell-1,
\end{equation}

or $B_1 = B_\ell$ and 
\begin{equation}\label{bineq2}
  s\geq 2+\ell(  b_1-1).
\end{equation}

Recall that by the maximality of $V(P)$, all edges of $H'$ containing $u_1$ or $u_\ell$ are contained in $V(Q) \cup V(P)$. For $j \in \{1,\ell\}$, define $A_j = N_{H'}(u_j) \cap V(Q)$ and $a_j = |A_j|$. By Claim~\ref{noconsecutive}, $A_j$ contains no consecutive vertices of $Q$.

\medskip
{\bf Case 1}: $A_1 = \emptyset$. Then all edges in $H'$ containing $u_1$ are contained in $V(P)$.

\medskip
{\bf Case 1.1}: $r=t'$. Since $\ell \leq t'-1$, no edge can be contained entirely in $V(P)$. Thus $u_1$ must only be contained in edges of $Q$ and $P$.

Then $b_1 \geq \delta(H)-|E(P)| = 1 + {t' \choose r-1} - (\ell -1) = t'-\ell +2$. If $(\ref{bineq1})$ holds, then 
\[n-\ell\geq s\geq 2(t'-\ell+2)+\ell-1=2t'+3-\ell \geq   n+2-\ell,\] 
a contradiction.

If instead $(\ref{bineq2})$ holds, then 
\begin{equation}\label{b1ineq2}
n\geq \ell+s\geq \ell+2+\ell((t'-\ell+2)-1)=2+\ell(t'-\ell+2).  
\end{equation}
Since
 $2 \leq \ell \leq t'-1$,    $\ell(t'-\ell+2)\geq 2(t'-2+2)\geq n-1$,    contradicting~\eqref{b1ineq2}.
 
\medskip
{\bf Case 1.2}: $3 \leq r \leq t'-1$. The number of edges in $H'$ containing $u_1$ and contained in $V(P)$ is at most ${\ell-1 \choose r-1}$. Thus, 
\[b_1 \geq \delta(H) - {\ell-1 \choose r-1} - |E(P)| = 1 + {t' \choose r-1} - {\ell-1 \choose r-1} - (\ell-1)\]\[ \geq  1 + {t' \choose 2} - {\ell-1 \choose 2} - (\ell-1) = \frac{(t'+\ell-2)(t'-\ell+1)}{2} - \ell +2.
\]

If (\ref{bineq1}) holds, then 
\[
\frac{(t'+\ell-2)(t'-\ell+1)}{2} - \ell +2 \leq \frac{s-\ell+2}{2} \leq \frac{n}{2}-\ell+1.
\]

However, $\frac{(t'+\ell-2)(t'-\ell+1)}{2} \leq t' -1$ implies that $0 \geq t'^2-\ell^2-3t'+3\ell = (t'-\ell)(t'+\ell-3)$. This cannot hold because $2 \leq \ell \leq t'-1$ and $t' \geq 3$, so $\frac{(t'+\ell-2)(t'-\ell+1)}{2} \geq t' > n/2-1,$ a contradiction.

If instead (\ref{bineq2}) holds, then 

\[
\frac{(t'+\ell-2)(t'-\ell+1)}{2} - \ell +2 \leq \frac{s-2}{\ell} +1 \leq \frac{n-2}{\ell}.
\]

However, we have $\ell + \frac{n-2}{\ell} \leq \frac{n+\ell^2-2}{2} \leq \frac{n}{2}+1 $, and thus $\frac{n-2}{\ell} \leq \frac{n}{2}-\ell+1 < \frac{(t'+\ell-2)(t'-\ell+1)}{2} - \ell +2$.

{\bf Case 2}: $A_1 \neq \emptyset, B_\ell \neq \emptyset$. 
Let $B = \bigcup_{e_j \in B_\ell} \{v_j, v_{j+1}\}$.
By Claim~\ref{distance3}, Lemma~\ref{indep2} and the facts $s \leq n-\ell \leq 2t'+1 - \ell$, $B_\ell \neq \emptyset$, $|B| \geq b_\ell + 1$, and $B \cap A_1 = \emptyset$, we have 

\begin{equation}\label{a1}
a_1 \leq \lfloor \frac{(2t'+1 - \ell) - b_\ell-1 - \ell + 1}{2} \rfloor = t'-\ell + 1 - \lceil b_\ell /2\rceil.\end{equation}

Recall that we assumed  $b_1 \leq b_\ell$. Therefore \begin{equation}
\label{u1deg}
d(u_1) \leq {a_1 + |V(P) - \{u_1\}| \choose r-1} + b_1 + |E(P)| \leq {a_1 + \ell-1 \choose r-1} + b_\ell + \ell-1,
\end{equation}
with equality only if $u_1$ belongs to every edge of $P$, and $b_1 = b_\ell$. 

Combining~\eqref{u1deg} and~\eqref{a1}, we obtain
\begin{equation}\label{u1deg2}{t' \choose r-1} + 1 \leq d(u_1) \leq {t' - \lceil b_\ell/2 \rceil \choose r-1} + b_\ell + \ell -1.\end{equation}

{\bf Case 2.1}: $r = t'$. Since $A_1 \neq \emptyset$, we need $a_1 \geq r - |V(P)| = t' - \ell$. By~\eqref{a1}, $1\leq b_\ell\leq 2$, and $a_1 + \ell-1 =r-1$. Then from~\eqref{u1deg2}, we get

\[t' + 1 = {t' \choose r-1} + 1 \leq d(u_1) \leq b_\ell + \ell \leq 2 + \ell \leq 2 + (t'-1).\]
This gives a contradiction unless $b_1 = b_\ell = 2$ and $\ell = t'-1$ (so $s = |V(Q)| \leq t'+2 = \ell+3$). But then there is no way to fit two edges in $B_1$ and two edges in $B_\ell$ without violating Claim~\ref{distance}.

{\bf Case 2.2}: $3\leq r \leq t'-1$. If $t'-\lceil b_\ell/2 \rceil \leq r-1$, then as in the previous subcase,
$d(u_1) \leq b_\ell + \ell.$ Since $A_1 \neq \emptyset$ and $\ell \geq 2$, in order not to violate Claim~\ref{distance3} we need $b_\ell \leq |E(Q)| - 2 = s-3$.  Therefore $d(u_1) \leq s-3 + \ell \leq n-3 \leq 2t'-2$.
 When $t' \geq 4$ ($n \geq 8$), we have ${t' \choose r-1} +  1 > 2t'-2$, a contradiction. In the remaining cases $6\leq n \leq 7$, $r\geq 3$ implies $r=3 = t'$ which was handled in the previous subcase.

So we may assume $t' - \lceil b_\ell/2 \rceil >r-1$ and therefore from~\eqref{u1deg2} we get, 
\begin{equation}\label{u1deg3}{t' \choose r-1} - {t' - \lceil b_\ell/2 \rceil \choose r-1} = {t'-1 \choose r-2} + \ldots + {t'-\lceil b_\ell/2 \rceil \choose r-2} \leq b_\ell+ \ell-2 \leq b_\ell + t'-3.\end{equation}

Here we use the fact that $\ell \leq t'-1$ and ${t' \choose r-1} - {t'- c \choose r-1} = {t'-1 \choose r-2} + {t'-2 \choose r-2} + \ldots + {t'-c \choose r-2}$ for any positive integer $c\leq \lceil b_1/2\rceil$.

Let $f(x) = {t' -1 \choose r-2} + \ldots + {t'-\lceil x/2 \rceil \choose r-2}$ and $g(x) = x + t' - 3$. 
For $x \in \{1,2\}$ and $r \geq 3$, $f(x) = {t'-1 \choose r-2} \geq t'-1 \geq g(x)$ with equality only if $r = 3$ and $x = 2$. For integers $2 < x \leq b_\ell$, $g(x) = g(2) + (x-2)\leq g(2) + 2\lceil (x-2)/2 \rceil$, and 
$$f(x) \geq f(2) + {t'-2 \choose r-2} + \ldots + {t'- \lceil x/2 \rceil \choose r-2}.$$
Each of these terms is at least 2, so $f(x) \geq f(2) + 2\lceil x/2 \rceil$. So $f(b_\ell) > g(b_\ell)$ if $b_\ell \neq 3$, contradicting~\eqref{u1deg3}. The final case is $b_\ell=2$. Moreover, we also get a contradiction if equalities in~\eqref{u1deg}--\eqref{u1deg3} do not hold. In this case, we must have $b_1 = b_\ell = 2$, and $\ell = t'-1$ (so $s \leq t'+2$). But then there is no way to fit two edges in $B_1$ and two edges in $B_\ell$ without violating Claim~\ref{distance}.
%
%
%
%
%
%
This finishes Case 2.

{\bf Case 3}: $B_1, B_\ell = \emptyset$.
\bigskip
Let us show that
\begin{equation}\label{a2}
\mbox{\em if  $B_1, B_\ell = \emptyset$, then $a_1,a_\ell \geq t'-\ell+1\geq 2$.}
\end{equation}
Indeed, if $a_1\leq t'-\ell$, then $d_{H'}(u_1)\leq {t'-\ell+(|V(P)|-1)\choose r-1}={t'-1\choose r-1}$. So, since $B_1=\emptyset$ and $\ell\leq t'-1$,
$$d_H(u_1)= d_{H'}(u_1)+|E(P)|\leq {t'-1\choose r-1}+\ell-1\leq {t'\choose r-1},$$
a contradiction. The same argument works if $a_\ell\leq t'-\ell$.

Similarly, if $B_1, B_\ell = \emptyset$ and $i\in \{1,\ell\}$, then at least two edges in $H'$ containing $u_i$ are not subsets of $V(P)$. Indeed, otherwise
$$d_H(u_1)= d_{H'}(u_1)+d_{H_P}(u_1)\leq 1+{\ell-1\choose r-1}+\ell-1< {t'\choose r-1}.$$

Let $f, f'$ be distinct edges of $H'$ such that for distinct $i, j,$ we have $\{v_i, u_1\} \subset f, \{v_j, u_\ell\} \subset f'$.  By Claim~\ref{distance2}, $|j - i| \geq \ell+1$, and by Claim~\ref{noconsecutive}, $A_1$ contains no consecutive vertices in $Q$. Without loss of generality, there exists $v_i \in A_1$, $v_j \in A_\ell$ with $j \geq i+\ell+1$. 
 If $a_1\geq t'-\ell+2$, then
 $$n-\ell\geq s \geq |A_1|+|\{v_{k+1}: v_k \in A_1, k\neq i, s\}|+|\{v_{i+1}, \ldots, v_{i+\ell}\}|$$ $$ \geq 2(t'-\ell+2)-2+\ell\geq 2\frac{n-1}{2}-\ell+2=n-\ell+1,$$
 a contradiction. Hence by~\eqref{a2},  
 \begin{equation}\label{a=t'} 
 a_1 = a_\ell=t'-\ell+1 \geq 2.
\end{equation}

We also prove that
\begin{equation}\label{a2''}
\parbox{14cm}{\em if  $B_1, B_\ell = \emptyset$, then  every $v_i \in A_1$  is contained in at least two common edges of $H'$ with $u_1$, and similar with $u_\ell$.}
\end{equation}
Indeed, otherwise by~\eqref{a=t'},

\begin{eqnarray*}d_H(u_1) &=& d_{H'}(u_1)+d_{H_P}(u_1)\\
&\leq & {|A_1 \cup V(P) -\{u_1\}| \choose r-1} - |\{f \subseteq A_1 \cup V(P): v_i\in f, |f| = r-1\}| + 1 +|E(P)|\\
& \leq & {t' -\ell+1 + (\ell-1) \choose r-1} - {t' -\ell+1 + (\ell-1) -1 \choose r-2} + 1 + \ell-1\\
& \leq & {t' \choose r-1} -(t'-1) +1+ \ell-1 \\
&\leq  & {t' \choose r-1} < \delta(H).
\end{eqnarray*}

This implies that for each $v_i\in A_1$ and $v_j \in A_\ell$, there exist distinct edges $f, f' \in H'$ such that $\{u_1, v_i\} \in f, \{u_\ell, v_j\} \in f'$, and hence $|j-i| \geq \ell+1$ by Claim~\ref{distance2}.

\smallskip
{\bf Case 3.1}: $B_1, B_\ell = \emptyset$ and $A_1 \neq A_\ell$. 
%
Without loss of generality, $|A_\ell| \geq |A_1|$. 
Since $A_1, A_\ell$ are independent sets and $A_\ell - A_1$ is nonempty, by Lemma~\ref{verc2}, $s \geq 2a_1 + 2|A_\ell-A_1| + (\ell+1)-3$. Hence
\begin{equation}\label{u1deg4} a_1 \leq \lfloor \frac{s - 2|A_\ell - A_1| - \ell +2}{2} \rfloor \leq  \lfloor\frac{(2t'+1-\ell) - 2|A_\ell - A_1| - \ell +2}{2} \rfloor = t'-\ell -|A_\ell - A_1|+1.
\end{equation}

We have
\[d(u_1) \leq {a_1 + \ell-1 \choose r-1} + |E(P)| \leq {t' - |A_\ell - A_1| \choose r-1} + \ell-1.\]
Since $|A_\ell - A_1|\geq 1$ and $\ell \leq t'-1$, this quantity is strictly less than ${t' \choose r-1} + 1$, a contradiction.

\bigskip
{\bf Case 3.2}: $B_1, B_\ell = \emptyset$ and $A_1 = A_\ell$. 
If $H'$ has no edges containing both, $u_1$ and $u_\ell$, then by the case and~\eqref{a=t'}, 
$$d_H(u_1)=d_{H'}(u_1)+d_{H_P}(u_1)\leq {t'\choose r-1}-|\{f\subseteq A_1 \cup V(P): |f| = r-1, u_\ell \in f\}| + |E(P)|$$ $$ ={t' \choose r-1} - {t'-1\choose r-2}+(\ell-1)\leq {t'\choose r-1},
$$
a contradiction. So suppose there is $f_0\in E(H')$ containing $\{u_1,u_\ell\}$. Let $$P_j = u_j, f_{j-1}, \ldots, f_1, u_1, f_0, u_\ell, f_{\ell-1}, \ldots, u_{j+1}$$ denote the path obtained from $P$ by adding $f_0$ and deleting $f_j$. By definition, for each $1\leq j\leq \ell-1$, the pair $(Q,P_j)$ is also a best pair. This yields that each $f_j$ is contained in $V(Q) \cup V(P)$. Moreover, 
for each such $j$ we have Case 3.2. By~\eqref{a2''}, deleting $f_0$ from $H'$ does not change $A_1$. It follows that
$A_1=A_2=\ldots=A_\ell$, and hence each $f_j$ is contained in $A_1\cup V(P)$. So by~\eqref{a=t'},
$d_H(u_1)\leq {(t'-\ell+1)+(\ell-1)\choose r-1}$, a contradiction.
\end{proof}

\section{Proof of Theorem~\ref{main} for $r > n/2$}

In this section, we complete the proof of Theorem~\ref{main} by showing that if $r > n/2 \geq 3$ and $\delta(H) \geq r-1$
or $r=3$, $n=5$ and $\delta(H) \geq 3$, then $H$ is hamiltonian-connected.

\begin{proof}[Proof of Theorem~\ref{main} for $r > n/2$.]
Suppose that  an $r$-graph $H$ with $\delta(H)\geq r-1$ has no hamiltonian  $x,y$-path for some
  $x,y\in V(G)$.  Let  $(Q, P)$ be a best  pair of two vertex-disjoint paths $Q$ and $P$ such that $Q$ is a $x, y$-path. 
  
  It is straightforward to check that the theorem is satisfied when $n=4$, $r=3$, $\delta(H) = 2$, so we may assume $\delta(H) \geq 3$.
  
  Since by Lemma~\ref{r+1path}, $s\geq r+1$ and
  $r\geq \left\lceil\frac{n+1}{2}\right\rceil$,  we have $\ell\leq n-s\leq \left\lfloor\frac{n-3}{2}\right\rfloor$ and
  \begin{equation}\label{r-l}
  r-\ell\geq \left\lceil\frac{n+1}{2}\right\rceil-\left\lfloor\frac{n-3}{2}\right\rfloor\geq 2.
  \end{equation}
  
  \medskip
  {\bf Case 1:} $\ell=1$. As in Section~\ref{step2}, in this case every edge $g \in H'$ contains at most one vertex outside of $V(Q)$.

  \medskip
  {\bf Case 1.1:} There are two edges $g,g'\in E(H')$ containing $u_1$. Then $|(g\cup g')\cap V(Q)|\geq r$, and no two vertices of $g\cup g'$ are consecutive on $Q$. It follows that $s\geq 2r-1\geq 2\frac{n+1}{2}-1=n$, a contradiction to $s\leq n-\ell$.

  \medskip
  {\bf Case 1.2:} There is exactly one edge $g\in E(H')$ containing $u_1$. Since  $\delta(H)\geq 3$, at least two edges of $H_Q$ contain $u_1$. By Claim~\ref{notneighbor}, $g$ does not intersect the sets $\{v_i,v_{i+1}\}$ such that
  $u_1\in e_i$. On the other hand, since no two vertices in $g$ are consecutive on $Q$, the $r-1$ vertices of $g\cap V(Q)$ intersect at least $2(r-1)-2\geq n-3$ sets $\{v_i,v_{i+1}\}$. This contradicts the fact that $Q$ has  $s-1\leq n-2$ pairs
  $\{v_i,v_{i+1}\}$.

  \medskip
  {\bf Case 1.3:} All edges  containing $u_1$ are in $B_1$, and some edge  $g\in H'$ is contained in $V(Q)$. Then $d(u_1) = b_1 \geq r-1$.
%
%
By Claim~\ref{forbidden2}, $r-1\leq b_1 \leq s-r + 1$ and therefore $n\leq 2r-1 \leq s+1 \leq n$. This implies that the ``equality" part of Claim~\ref{forbidden2}(i) holds, and so $g =\{v_1, \ldots, v_i\} \cup \{v_{j+1}, \ldots, v_s\}$ and $B_1 = \{e_i, \ldots, e_j\}$ for some $i<j$. 
  In particular, by symmetry we may assume that $i>1$. Let $Q'$ be the path obtained by replacing $e_1$ with $g$. We get a new best pair $(Q',P)$ with $e_1$ playing the old role of $g$. As $B_1$ does not change, Claim~\ref{forbidden2} asserts $g =e_1$, a contradiction.

  \medskip
  {\bf Case 1.4:} All edges  containing $u_1$ are in $B_1$, and no edges in  $ H'$ are contained in $V(Q)$.
  Again, $|B_1|\geq r-1$.  Since $|E(H)|\geq n-1$, there is an edge $g\in E(H')$. Since $\ell=1$ and Case 1.3 does not hold,
  $|g\cap V(Q)|=r-1$. Then $g$ has a vertex $w$ outside of $V(Q)\cup \{u_1\}$, so $s\leq n-2$.
  
If there is another edge $g'\in E(H')$ containing $w$, then there could not be consecutive vertices $v_i,v_{i+1}$ in $Q$
such that one of them is in $g$ and the other in $g'$.
Hence
 the sets $A=g\setminus \{w\}$ and 
  $B=g'\setminus \{w\}$ satisfy condition~\eqref{dis} for $q=2$ in Lemma~\ref{ver-new}. Since $g'\not\subseteq g$
 and $ g\not\subseteq g'$, 
  Lemma~\ref{ver-new}(ii) for $q=2$ yields $s\geq |A|+|B|+q-1\geq 2r-1\geq n$. This contradicts the fact that $s\leq n-2$.
  
  Otherwise, $w$ belongs to some $r-2$ edges $e_{i_1},\ldots,e_{i_{r-2}}$. Let $A=\bigcup_{j=1}^{r-2}\{v_{i_j},v_{i_j+1}\}$.
  Then $|A|\geq r-1$.
  By Claim~\ref{notneighbor}, $g\cap A=\emptyset$.  Hence $s\geq (r-1)+(r-1)\geq 2\frac{n+1}{2}-2=n-1$, contradicting $s\leq n-2$.

  \medskip
  {\bf Case 2:} $2\leq\ell\leq \left\lfloor\frac{n}{2}\right\rfloor-1$.
  
 \medskip
  {\bf Case 2.1:} $a_1\geq 1$ and $b_\ell\geq 1$. Let $g\in E(H')$ contain $u_1$. Then $g\subset V(Q)\cup V(P)$ and $|g\cap V(Q)|\geq r-\ell$. Since $\ell\geq 2$, by Lemma~\ref{indep2} with $q=\ell$,
  $$r-\ell\leq \frac{s-1-\ell+1}{2}\leq \frac{n}{2}-\ell,
  $$
  contradicting $r>\frac{n}{2}$.
  
  \medskip

If $b_1 \leq 1$, then $u_1$ is contained in at least one edge in $H'$ and so $a_1 \geq r - \ell \geq 2$ by~\eqref{r-l}. Either way, $a_1 + b_1 \geq 2$ and similarly $a_\ell + b_\ell \geq 2$. 
  By symmetry, the following two subcases remain.
  
  {\bf Case 2.2:} $b_1=b_\ell=0$. Then $d_{H'}(u_1)\geq \delta(H) - |E(P)|\geq (r-1)-(\ell-1)\geq 2$ and similarly 
   $d_{H'}(u_\ell)\geq  2$. Let $g_1\in E(H')$ contain $u_1$ and $g_\ell\in E(H')-g_1$ contain $u_\ell$.
 Let $A=g_1\cap V(Q)$  and $B=g_\ell\cap V(Q)$. Then $A$ and $B$ satisfy condition~\eqref{dis} for $q=1+\ell$ in Lemma~\ref{ver-new}.
  
 Also, $|A|\geq r-\ell$ with equality only if $g_1\supset V(P)$, and the same holds for $B$. If $|A|=|B|=r-\ell$, then $A\neq B$ because $g_1\neq g_\ell$. In this case, by Lemma~\ref{ver-new}(ii) for $q=1+\ell$,
 $s\geq 2(r-\ell)+(1+\ell)-1=2r-\ell$. Since $s\leq n-\ell$, this contradicts $r>n/2$. Similarly, if $\max\{|A|,|B|\}\geq r-\ell+1$ and $A\neq B$, then  by Lemma~\ref{ver-new}(ii) for $q=1+\ell$,
 $s\geq (r-\ell)+(r-\ell+1)+(1+\ell)-2=2r-\ell$. So, we get the same contradiction. 
 
 Finally, suppose $A=B$. Since $g_1\neq g_\ell$, this implies $|A|\geq r-\ell+1$. Hence 
 by Lemma~\ref{ver-new}(i) for $q=1+\ell$,  we get
 $$n-\ell\geq s \geq 1+(\ell+1)(r-\ell),$$
 which yields
 \begin{equation}\label{Case22}
 \ell(r-\ell)\leq n-r-1.
 \end{equation}
 Since $2\leq \ell\leq r-2$, for fixed $n$ and $r$,  the LHS in~\eqref{Case22} is at least 
 $2(r-2)$. Thus~\eqref{Case22} implies $3r\leq n+3$. But $3r\geq (n+1)+r\geq n+4$, a contradiction.

  {\bf Case 2.3:} $a_1=a_\ell=0$.  Similarly to Case 2.2,  $b_1\geq r-\ell$ and $b_\ell\geq r-\ell$.
 Let $A=\{v_iv_{i+1}: u_1\in e_i\}$  and $B=\{v_jv_{j+1}: u_\ell\in e_j\}$.
  Then $A$ and $B$ satisfy condition~\eqref{dise} for $q=\ell$ in Lemma~\ref{ed-new}.
 
 Also, $|A|\geq r-\ell$ with equality only if $u_1\in f_j$ for all $j$, and the same holds for $B$ (with $u_\ell$ in place of $u_1$). If $|A|=|B|=r-\ell$ and $A\neq B$, then by Lemma~\ref{ed-new}(ii) for $q=\ell$,
 $s-1\geq 2(r-\ell)+\ell-1=2r-\ell-1$. Since $s\leq n-\ell$, this contradicts $r>n/2$. Similarly, if $\max\{|A|,|B|\}\geq r-\ell+1$ and $A\neq B$, then  by Lemma~\ref{ed-new}(ii) for $q=\ell$,
 $s-1\geq (r-\ell)+(r-\ell+1)+\ell-2=2r-\ell-1$. So, we get the same contradiction. 
 
 Finally, suppose $A=B$. Let  $B'=\bigcup_{\{j: u_1 \in e_j\}}\{v_j,v_{j+1}\}$. Since $A=B$, $|B'|\geq 2(r-\ell)$.
 Let $A'=f_1\cap V(Q)$. If $A'=\emptyset$, then $|V(H)-V(Q)-V(P)|\geq |f_1-V(P)|\geq r-\ell\geq 2$.
Then by Lemma~\ref{ed-new}(i) for $q=\ell$,  we get 
 $n-\ell-2\geq s\geq 2+\ell(r-\ell-1)$, which yields $\ell(r-\ell)\leq n-4$. Since $2\leq \ell\leq r-2$,
 the LHS of this inequality is at least $2(r-2)\geq n-3$,  a contradiction. 
  Thus $A'\neq \emptyset$. If $v_{i_1}\in A'\cap B'$,
say   $e_{i_1}\in B$, then we can replace edge $e_{i_1}$ in $Q$ by the path $v_{i_1},f_1,u_1,e_{i_1},v_{i_1+1}$, contradicting the choice of $(Q,P)$. Thus $A'\cap B'=\emptyset$. Moreover, similarly if
$i_1<i_2\leq i_1+\ell-2$, $v_{i_1}\in A'$ and $e_{i_2}\in B'$, then we can 
replace  the subpath $v_{i_1},e_{i_1},v_{i_1+1},\ldots,v_{i_2}$ of $Q$ with the longest path
 $v_{i_1},f_1,u_2,f_2,u_3,\ldots,u_\ell, e_{i_2},v_{i_2+1}$, a contradiction again. It follows that
 $s\geq |A'|+|B'|+\ell-2$. If $|A'|=1$, then $s\leq n-\ell-1$, and therefore
 $n-\ell\geq |B'|+\ell\geq 2r-\ell >n-\ell$, a contradiction. Otherwise if $|A'| \geq 2$, then $s \geq |B'| + \ell \geq 2r-\ell > n-\ell$ again.

\end{proof}
 
 \section{Concluding remarks}
 

1. A number of theorems on graphs, in particular, Theorem~\ref{o2}, give sufficient conditions for the existence 
of   hamiltonian cycles in terms of $\sigma_2(G)=\min_{uv\notin E(G)}d(u)+d(v)$. Partially, this is because many proofs of bounds in terms of the
 minimum degree also work for $\sigma_2(G)$. It seems this is not the case for $r$-graphs when $r\geq3$. Moreover the degree of a vertex in an $r$-graph can be interpreted in different ways: the number of edges containing the vertex or the number of vertices in its neighborhood. Defining a suitable analog of $\sigma_2(G)$ for hypergraphs is unclear.  
  For example, if $n=2r$, then there are $n$-vertex $r$-graphs with $6$ edges in which every two vertices are in a common edge (e.g., a blow up of a $K_4$), so counting the sizes of the neighborhoods is not a useful parameter at least for large $r$. On the other hand for small $r$, the hypergraph consisting of a  $K_{n-1}^r$ and one additional edge satisfies $d(u) + d(v) \geq {n-2 \choose r-1} + 1$ for every pair of vertices and is not hamiltonian. While it is likely possible to prove an Ore-type theorem using this bound, this quantity is significantly larger than the sufficient minimum degree condition $\delta(H) \geq {\lfloor (n-1)/2 \rfloor \choose r-1}+1$  needed for hamiltonicity, and so such a result may not be very meaningful. 
It would be interesting to find some analog of $\sigma_2(G)$ for $r$-graphs that is both natural and nontrivial for a given range of $r$.

\medskip
2. Given $k\geq 2$, a (hyper)graph $G$ is {\em $k$-path-connected} if for any distinct $x,y\in V(G)$, there is an
$x,y$-path with at least $k$ vertices. In these terms, an $n$-vertex (hyper)graph is hamiltonian-connected exactly when it is $n$-path-connected. It would be interesting to find exact restrictions on the minimum degree of an $n$-vertex $r$-graph $G$ providing that $G$ is $k$-path-connected for $r<k<n/2$.

\medskip
3. Call a graph $G$ {\em $1$-extendable} if for each edge $e\in E(G)$,   $G$ has a hamiltonian cycle containing $e$. Thus Theorem~\ref{o2} yields that for $n\geq 3$ each $n$-vertex  graph $G$  with $\delta(G)\geq (n+1)/2$ is $1$-extendable.  Also, one can define
{ $1$-extendable} hypergraphs in several ways. One natural definition would be: 
An $r$-graph $G$ {\em $1$-extendable} if for each edge $e\in E(G)$ and any two vertices $u,w\in e$,   $G$ has a hamiltonian cycle $C=v_1,e_1,v_2,\ldots,v_n,e_n,v_1$ such that $e_1=e$, $v_1=u$ and $v_2=w$.

For $r=2$, this definition coincides with the original definition of $1$-extendable graphs, but for $r\geq 3$ the claim that 
each hamiltonian-connected $r$-graph is $1$-extendable is not true: as we have seen in 
Section~\ref{hcnh}, hamiltonian-connected $r$-graphs do not need to be even just hamiltonian.
On the other hand, trivially if each $n$-vertex $r$-graph with minimum degree at least $d$ is hamiltonian-connected, then each $n$-vertex $r$-graph with minimum degree at least $d+1$ is $1$-extendable. So, Theorem~\ref{main} yields
the following.

\begin{cor}\label{1ext} Let $n \geq r \geq 3$. Suppose $H$ is an $n$-vertex, $r$-graph such that (1) $r \leq n/2$ and $\delta(H) \geq {\lfloor n/2 \rfloor \choose r-1} + 2$, or (2) $r > n/2\geq 3$  and $\delta(H) \geq r$, or (3) $r=3$, $n=5$ and $\delta(H) \geq 4$.
Then $H$ is  $1$-extendable. 
\end{cor}

When $r>n/2\geq 3$, the bound in Corollary~\ref{1ext} is exact, but when $3\leq r\leq n/2$ or $r=3$ and $n=5$ it probably can be improved by $1$.

\medskip
4. 
P\'osa~\cite{Posa} considered the following generalization of $1$-extendable graphs. Given a linear forest 
 (i.e., a set of vertex-disjoint paths)
 $L$,  call a graph $G$   {\em $L$-extendable} if   $G\cup L $ has a hamiltonian cycle containing 
 all edges of $L$.  P\' osa~\cite{Posa}   proved that for each $n>\ell\geq 0$ and every linear forest $L$ with $\ell$ edges, each $n$-vertex graph $G$ with $\sigma_2(G)\geq n+\ell$ is $L$-extendable. This is a far reaching generalization of Theorem~\ref{o2}.
 It also implies that if~$\sigma_2(G)\geq n+\ell$, then $G$ is $\ell+1$-hamiltonian-connected (take $L$ to be a path on at most $\ell$ vertices). 
 
One may consider different hypergraph definitions of being $L$-extendable for a given graph linear forest $L$.  
For example, given a positive integer $\ell$,
 we can say that a hypergraph $G$ is
$\ell$-{\em extendable} if for every choice of $\ell+1$ vertices $u_1,\ldots,u_{\ell+1}$ and $\ell$ edges $g_1,\ldots,g_\ell$ in $G$ such that $\{u_i,u_{i+1}\}\subset g_i$ for all $i\in [\ell]$, $G$ has a hamiltonian cycle $C=v_1,e_1,v_2,\ldots,v_n,e_n,v_1$ such that  $v_i=u_i$ for all $i\in [\ell+1]$ and $e_j=g_j$ for all $j\in [\ell]$.
Exact bounds on the minimum degree providing that an $n$-vertex $r$-graph is $\ell$-{extendable} seem difficult for general $\ell$ but probably are feasible for very small or very large $\ell$.

\medskip
5. Chartrand,  Kapoor and Lick~\cite{CKL} proved analogs of Dirac's Theorem and its generalizations by Ore~\cite{ore} and P\' osa~\cite{Posa} for $\alpha$-{\em hamiltonian graphs}, that is, the graphs that are hamiltonian after deleting any set of at most $\alpha$ vertices. If in this definition we replace ``at most" with ``exactly", the class of the graphs satisfying the definition may change. For example, after deleting any vertex from Petersen Graph, the remaining graph is hamiltonian. Lick~\cite{Lick} proved similar exact results for $\alpha$-{\em hamiltonian-connected graphs}, that is, the graphs that are hamiltonian-connected after deleting any set of at most $\alpha$ vertices.
The ideas and tricks in~\cite{KLM} and this paper may be used to try to find exact or close to exact  bounds
on minimum degree in an $n$-vertex $r$-graph $G$ ensuring that $G$ is $\alpha$-hamiltonian/$\alpha$-hamiltonian-connected. As with graphs, the answers for the definitions with ``at most" and ``exactly" may differ.

\bigskip
{\bf Acknowledgment.} We thank the anonymous referees for their helpful comments.


\begin{thebibliography}{99}
\small


\bibitem{BGHS}
J.-C. Bermond, A. Germa, M.-C. Heydemann, D. Sotteau: Hypergraphes Hamiltoniens,
in \emph{Probl\`emes combinatoires et th\`eorie des graphes} (Colloq. Internat.
CNRS, Univ. Orsay, Orsay, 1976). Colloq. Internat. CNRS,  \textbf{260} (CNRS,
Paris, 1978), 39--43.


\bibitem{CKL}
G. Chartrand, S. F. Kapoor, D. R. Lick:
$n$-Hamiltonian graphs,
{\em J. Combinatorial Theory} {\bf 9} (1970), 308--312.


\bibitem{CEP}
D. Clemens, J. Ehrenm\" uller, and Y. Person, A Dirac-type theorem for Hamilton
Berge cycles in random hypergraphs, {\em Electron. J. Combin.} \textbf{27} (2020), no. 3, Paper No. 3.39, 23 pp.


\bibitem{CP}
M. Coulson, G. Perarnau:
A Rainbow Dirac's Theorem, 
{\em SIAM J. Discrete Math.} \textbf{34} (2020), no. 3, 1670--1692.

%
\bibitem{dirac}
G. A. Dirac:
Some theorems on abstract graphs, \emph{Proc. London Math. Soc.} {\bf 2}(3)
(1952), 69--81.
%
%
%
%
%
%
%
%


\bibitem{KLM}
A. Kostochka, R. Luo, G. McCourt:
Dirac's Theorem for hamiltonian Berge cycles in uniform hypergraphs, 
{\tt  arXiv:2109.12637}, (2022), 23 pp.
%

\bibitem{Lick}
D. R. Lick: A sufficient condition for Hamiltonian connectedness. {\em J. Combinatorial Theory}  {\bf 8} (1970), 444--445. 

\bibitem{MHG}
Y. Ma, X. Hou, J. Gao:
A Dirac-type theorem for uniform hypergraphs, 
{\tt  arXiv:2004.05073}, (2020), 17 pp.


\bibitem{ore}
O. Ore:
Note on Hamilton circuits, {\em 	Amer. Math. Monthly} \textbf{67}(1) (1960), 55.

\bibitem{ore2}
O. Ore:  Hamilton connected graphs. {\em J. Math. Pures Appl.} {\bf 42}(9) (1963), 21--27. 

\bibitem{Posa}   L. P\' osa:  A theorem concerning Hamilton lines, {\em Magyar Tud. Akad. Mat. Kutat\'o Int. K\"ozl.} {\bf 7} (1962), 225--226.


\bibitem{SN} N. Salia, P\' osa-type results for Berge hypergraphs,  {\tt arXiv:2111.06710v2},  15 p, 2021.

\end{thebibliography}
\end{document}